\documentclass{amsart}

\usepackage{amsmath,amsfonts,amssymb,amsthm,graphicx}

\newtheorem{theorem}{Theorem}[section]
\newtheorem*{theorem*}{Theorem}
\newtheorem{lemma}[theorem]{Lemma}
\newtheorem*{lemma*}{Lemma}
\newtheorem{corollary}[theorem]{Corollary}
\newtheorem{proposition}[theorem]{Proposition}
\newtheorem{claim}[theorem]{Claim}
\newtheorem*{bigthm}{Main Theorem}

\theoremstyle{remark}
\newtheorem{remark}[theorem]{Remark}
\newtheorem*{remark*}{Remark}
\newtheorem{example}[theorem]{Example}

\theoremstyle{definition}
\newtheorem{definition}[theorem]{Definition}

\def\C{\mathbb C}
\def\Z{\mathbb Z}
\def\Q{\mathbb Q}
\def\N{\mathbb N}
\def\zhs{\mathbb Z\text{HS}}
\def\qhs{\mathbb Q\text{HS}}
\def\xfn{(X_{f,n},0)}

\newcommand{\sing}[1]{(#1,0)}
\newcommand{\h}[1]{\widetilde{h_{#1}}}
\newcommand{\vbar}[1]{\overline{v_{#1}}}
\newcommand{\HH}[1]{h_{#1}\widetilde{h_{#1}}}
\newcommand{\ga}[1]{\Gamma_{A}(v_{#1})}
\newcommand{\gm}[1]{\Gamma_{-}(v_{#1})}
\newcommand{\dm}[1]{D_{-}(v_{#1})}
\newcommand{\da}[1]{D_{A}(v_{#1})}
\newcommand{\delm}[1]{\Delta_{-}(v_{#1})}
\newcommand{\dela}[1]{\Delta_{A}(v_{#1})}
\newcommand{\z}[2]{z_{#1,#2}}

\begin{document}

\title[topology of surface singularities $\{z^n=f(x,y)\}$]
{On the topology of surface singularities $\{z^n=f(x,y)\},$ for
$f$ irreducible}
\author[E.A. Sell]
{Elizabeth A. Sell}
\date{June 30, 2009}
\address{Department of Mathematics\\
Millersville University\\
P.O. Box 1002\\
Millersville, PA 17551-0302} \email{liz.sell@millersville.edu}

%\subjclass[2000]{14B05, 32S50} \keywords{surface singularity,
%suspension singularity, rational homology sphere, abelian cover,
%resolution graph, splice diagram, splice quotient}

\begin{abstract} The splice quotients are an interesting class of
normal surface singularities with rational homology sphere links,
defined by W. Neumann and J. Wahl. If $\Gamma$ is a tree of
rational curves that satisfies certain combinatorial conditions,
then there exist splice quotients with resolution graph $\Gamma.$
Suppose the equation $z^n=f(x,y)$ defines a surface $X_{f,n}$ with
an isolated singularity at the origin in $\mathbb{C}^3.$  For $f$
irreducible, we completely characterize, in terms of $n$ and a
variant of the Puiseux pairs of $f,$ those $X_{f,n}$ for which the
resolution graph satisfies the combinatorial conditions that are
necessary for splice quotients. This result is
topological;
whether or not $X_{f,n}$ is analytically isomorphic to a splice
quotient is treated separately. 
\end{abstract}

\maketitle

%%%%%%%%%%%%%%%%%%%%%%%%%%%%%%%%%%%%%%%%%%%%%%%%%%%%%%%%%%
%%%%%%%%%%%%%%%%%%%%% Section 1 %%%%%%%%%%%%%%%%%%%%%%%%%%%%%%%%
%%%%%%%%%%%%%%%%%%%%%%%%%%%%%%%%%%%%%%%%%%%%%%%%%%%%%%%%%%

\section{Introduction}\label{intro}

Let $\sing{X}\subset\sing{\C^k}$ be the germ of a complex analytic
normal surface singularity.  The intersection of $X$ with a
sufficiently small sphere centered at the origin in $\C^k$ is a
compact connected oriented three-manifold $\Sigma,$ called the
\textit{link} of $\sing{X},$ that does not
depend upon the embedding in $\C^k.$  
Let $\Gamma$ be the dual resolution graph of a good resolution of
the singularity. The
homeomorphism type of the link can be recovered from $\Gamma,$ and conversely, W.
Neumann proved that (aside from a few exceptions) the homeomorphism type of the link determines the minimal good resolution graph~\cite{N2}. 
One interesting class of normal surface singularities is the set
of those for which the link is a \textit{rational homology sphere}
($\qhs$) (i.e., $H_1(\Sigma,\Q)=0$). The link is a $\qhs$ if and
only if any good resolution graph $\Gamma$ of $\sing{X}$ is a tree of
rational curves.  

The work of Neumann and Wahl (described in $\S$\ref{n-w}; see also
\cite{nwuac} and \cite{GTE}) provides a method for generating
analytic data for singularities from topological data. Starting
with a resolution graph $\Gamma$ that satisfies certain
conditions, known as the ``semigroup and congruence conditions", 
one can produce defining equations for a normal surface
singularity with resolution graph $\Gamma.$ The singularities that
result from this algorithm are called \textit{splice quotients}.
If the link $\Sigma$ is a $\zhs$ ($H_1(\Sigma,\Z)=0$), then only
the semigroup conditions are relevant, and the singularities
produced by the algorithm are said to be \textit{of splice type}.
This work has led to a recent interest in the properties of splice
quotients and related topics (see \cite{Nem}, \cite{NemOk},
\cite{Okuma1}, \cite{Okuma2}, \cite{Stevens}), and there are still
many unanswered questions.

One of the first questions that arises is: How many singularities
with $\qhs$ link are splice quotients? There are two layers to the
problem - topological and analytic. If one has a singularity that
satisfies the necessary topological conditions (which depend only
on the resolution graph), then there exist splice quotients with
that topological type, but it is a separate issue to determine
whether the singularity is analytically isomorphic to a splice
quotient. 
Originally, one wondered whether all $\Q$-Gorenstein singularities
with $\qhs$ link would turn out to be splice quotients. However,
the first counterexamples were found in the paper of I.
Luengo-Velasco, A. Melle-Hern\'{a}ndez, and A. N\'{e}methi
\cite{Nem}. There, the authors give an example of a hypersurface
singularity for which the resolution graph does not satisfy the
semigroup conditions, and an example of a singularity for which
the semigroup and congruence conditions are satisfied, but the
analytic type is not a splice quotient. On the other hand, there
are nice classes of singularities for which \textit{all} analytic
types are splice quotients: weighted homogeneous singularities, as
shown by Neumann in \cite{N}, and rational and $\qhs$-link
minimally elliptic
singularities, as shown by T. Okuma in \cite{Okuma1}. 

A natural class of surface singularities to study after weighted
homogeneous, rational, and minimally elliptic is the class of
hypersurface singularities defined by an equation of the form
$z^n=f(x,y).$  If $\{f(x,y)=0\}$ defines a reduced curve with a
singularity at the origin in $\C^2,$ then for $n>1,$ the surface
$X_{f,n}:=\{z^n=f(x,y)\}$ has an isolated (hence normal)
singularity at the origin in $0\in\C^3.$ For $f$ irreducible, the
resolution graph of
$\xfn$ 
can be constructed from $n$ and a finite set of pairs of positive
integers associated to $f,$ known as the \textit{topological
pairs} $\{(p_i,a_i)~|~1\leq i\leq s\}$ defined in \cite{EN} (a
variant of the more commonly known \textit{Puiseux pairs}). The
topological pairs completely determine the topology of the plane
curve singularity. If there is only one topological pair ($s=1),$
then any such $\xfn$ with $\qhs$ link has the topological type of
a weighted homogeneous singularity, hence has the topological type
of a splice quotient. In \cite{Cass}, Neumann and Wahl prove that
the link of $\xfn$ is a $\zhs$ if and only if $f$ is irreducible
and all $p_i$ and $a_i$ are relatively prime to $n,$\footnote{In
\cite{Cass}, the result is incorrectly stated.  The pairs in
question are mistakenly identified as the Newton pairs instead of
the topological pairs.} and in that case, they prove in \cite{zhs}
that any such $\xfn$ is of splice type. That is, not only are the
semigroup conditions satisfied, but moreover, every $\xfn$ with
$\zhs$ link is isomorphic to one that results from Neumann and
Wahl's construction.

The main result of this paper is a complete characterization of
the $\xfn,$ with $f$ irreducible and $s\geq 2,$ that have a
resolution graph that satisfies the semigroup and congruence
conditions. 
For $f$ irreducible, there is an explicit criterion given by R.
Mendris and N\'{e}methi in \cite{MenNem}, in terms of $n$ and the
topological pairs, that determines when the link of $\xfn$ is a
$\qhs$ (see Proposition \ref{resgraph}). 
One can see that there
are plenty of $\xfn$ for which the link is a $\qhs$ but not a
$\zhs.$ From now on, whenever we are not referring to topological
pairs, the notation $(m,n)$ denotes the greatest common divisor of
the integers $m$ and $n.$ Our main result is the following

\begin{bigthm}Let $f$ be irreducible with
topological pairs $\{(p_i,a_i)~:~1\leq i\leq s\},$ with $s\geq2,$
and let $n$ be an integer greater than $1.$ Then $\xfn$ has $\qhs$
link and a good resolution graph that satisfies the semigroup and
congruence conditions if and only if either
\begin{itemize}\item[(i)] $(n,p_s)=1,~(n,p_i)=(n,a_i)=1$ for $1\leq
i\leq s-1,$ and $a_s/(n,a_s)$ is in the semigroup generated by
$\{a_{s-1},~p_1\cdots p_{s-1},~a_{j}p_{j+1}\cdots p_{s-1}~:~1\leq
j\leq s-2\},$ or \item[(ii)] $s=2,~p_2=2,~(n,p_2)=2,$ and
$(n,a_2)=(\frac{n}{2},p_1)=(\frac{n}{2},a_1)=1.$
\end{itemize}\end{bigthm}

It is somewhat surprising that so few $\xfn$ satisfy the
topological conditions, given the result in the $\zhs$ case. Aside from Case (ii), which is
rather restrictive, this result says that if any of the
topological pairs other than $a_s$ have factors in common with
$n,$ then $\xfn$ does not have the topological type of a splice
quotient. One could say that if $\xfn$ gets ``too far" from the
$\zhs$ case (for which all \textit{analytic} types are splice
quotients), it cannot even have the topology of a splice quotient.

If the resolution graph does satisfy the semigroup and congruence
conditions, a priori we do not know what the equations of the
splice quotients produced from the Neumann-Wahl algorithm look
like. Not only is it unclear whether or not $\xfn$ itself is a
splice quotient, but in fact, it is not even clear that there
exist splice quotients defined by any equation of the form
$z^n=g(x,y).$ It turns out that there do exist such splice
quotients; unfortunately, the length of the proof is such that it
cannot be included here. That result can be found in
\cite{paper2}. In the case of weighted homogeneous splice
quotients, it was shown in \cite{mythesis} that in general, not every
deformation with the same topological type is analytically
isomorphic to a splice quotient. Therefore, we expect that there
are few cases for which \textit{every} $\xfn$ of a given
topological type is a splice quotient.

Consider the following example.
\begin{example} Let $X_n:=\{z^{n}=y^5+(x^3+y^2)^2\}.$ The plane curve singularity
defined by $y^5+(x^3+y^2)^2=0$ is irreducible with two topological
pairs, $p_1=2,~a_1=3,~p_2=2,$ and $a_2=15.$ The link of
$\sing{X_n}$ is a $\qhs$ if and only if either $(n,2)=1$ or $(n,15)=1.$ We can say the following about $X_n:$
\begin{itemize}
\item If $n$ is relatively prime to $2,~3,$ and $5,$ then
$\sing{X_n}$ has $\zhs$ link and hence is of splice type. In fact,
we could replace $y^5+(x^3+y^2)^2$ by any curve with the same
topological pairs, and we would still have a singularity of splice
type. \item If $n$ is divisible by $3,$ the Main Theorem says that
$\sing{X_n}$ does not even have the topological type of a splice
quotient. \item If $n=5k,$ where $k$ is relatively prime to $2$
and $3,$ then $\sing{X_n}$ has the topology of a splice quotient
by Case (i) of the Main Theorem, and in fact,
$\sing{X_n}$ is itself a splice quotient \cite{paper2}. 
\item If $n=2k$, where $k$ is relatively prime
to $2,~3,$ and $5,$ then $\sing{X_n}$ has the topology of a splice
quotient
by Case (ii) of the Main Theorem. 
It is unclear whether or not $\sing{X_n}$ is a splice quotient.
However, if we replace $y^5+(x^3+y^2)^2$ by
$(x^3-y^2-y^3)^2-4y^5,$ which has the same topological pairs, it
is a splice quotient \cite{paper2}.
\end{itemize}\end{example}

The rest of this paper is entirely devoted to proving the Main
Theorem. In section \ref{n-w}, we provide a brief summary of the
work of Neumann and Wahl. Section \ref{top} contains a
description of the resolution graph and splice diagram for $\xfn.$
Some of the computations that are necessary for the proof of the
Main Theorem depend upon work done by Mendris and N\'{e}methi in
\cite{MenNem}; section \ref{resgraph_sec} is a reiteration of this
material. In section \ref{sgc_sec}, we analyze the semigroup
conditions for the splice diagram associated to $\xfn.$ Section
\ref{discgroup_sec} contains additional computations that are
needed for checking the congruence conditions.  Finally, in
section \ref{toptypes}, we use the computations from the previous
three sections to prove the Main Theorem.

\textit{Acknowledgements.} Special thanks are due to J. Wahl for
many very helpful conversations throughout the preparation of this
work.

%%%%%%%%%%%%%%%%%%%%%%%%%%%%%%%%%%%%%%%%%%%%%%%%%%%%%%%%%%
%%%%%%%%%%%%%%%%%%%%% Section 2 %%%%%%%%%%%%%%%%%%%%%%%%%%%%%%%%
%%%%%%%%%%%%%%%%%%%%%%%%%%%%%%%%%%%%%%%%%%%%%%%%%%%%%%%%%%

\section{The Neumann-Wahl algorithm}\label{n-w}

This section contains a summary of the method defined by Neumann
and Wahl in \cite{NW} to produce equations for the splice
quotients and their universal abelian covers; we refer to this
method as the \textit{Neumann-Wahl algorithm}. The algorithm
begins with a negative-definite graph $\Gamma$ that is a tree of
smooth rational curves (equivalently, the dual resolution graph
associated to a good resolution of a normal surface singularity
with $\Q$HS link)
and the \textit{splice diagram} $\Delta$ associated to $\Gamma.$ Splice
diagrams were introduced by Eisenbud and Neumann \cite{EN} for plane
curve singularities (building on work of Siebenmann), and later
generalized by Neumann and Wahl.
If $\Delta$ satisfies the ``semigroup conditions" (Definition
\ref{sgcdef}), then the algorithm produces a set of equations that
defines a family of isolated complete intersection surface
singularities. The algorithm also produces an action of the finite
abelian group $D(\Gamma),$ the discriminant group of $\Gamma,$ on the
coordinates used for the splice diagram equations. If $\Gamma$
satisfies further combinatorial conditions, the ``congruence
conditions" (Definition \ref{ccdef}), then one can choose a set of
splice diagram equations such that the discriminant group acts on every singularity $\sing{Y}$ in the family.
Furthermore, the quotient of $\sing{Y}$ by $D(\Gamma)$ is an isolated
normal surface singularity with resolution graph $\Gamma,$ and the
covering given by the quotient map is the universal abelian
covering (the maximal abelian covering that is unramified away
from the singular point). 

\begin{figure}[h]\centering
\includegraphics[width=\textwidth]{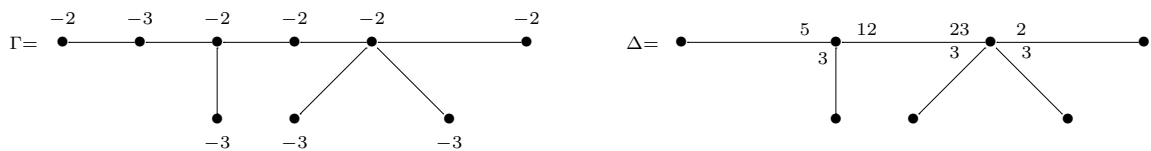}
\caption{A resolution graph $\Gamma$ and its associated splice diagram
$\Delta.$}\label{ex2_1}
\end{figure}

In a weighted graph, the \textit{valency} of a vertex is the
number of adjacent edges. A \textit{node} is a vertex of valency
at least three, a \textit{leaf} is a vertex of valency one, and a
\textit{string} is a connected subgraph that does not include a
node. The procedure for computing the splice diagram $\Delta$
associated to a resolution graph $\Gamma$ 
is as follows. First, omit the self-intersection numbers of
the vertices and contract all strings of valency two vertices in
$\Gamma.$ To each node $v$ in the resulting diagram $\Delta,$ we
attach a weight $d_{ve}$ in the direction of each adjacent edge
$e.$  
Remove the vertex in $\Gamma$ that corresponds to the node $v$ and the edge
that corresponds to $e,$ and let $\Gamma_{ve}$ be the remaining
connected subgraph that was connected to $v$ by $e.$ Then the
weight $d_{ve}=\det(-C_{ve}),$ where $C_{ve}$ is the intersection
matrix of the graph $\Gamma_{ve}.$  Figure \ref{ex2_1} contains a
simple example. Similarly, we define a subgraph $\Delta_{ve}$ of $\Delta$ as follows.
Remove $v$ and $e,$ and let $\Delta_{ve}$ be the remaining connected
subgraph that was connected to $v$ by $e.$ For any two vertices
$v$ and $w$ in $\Delta,$ the \textit{linking number} $\ell_{vw}$ is
the product of the weights adjacent to but not on the shortest
path from $v$ to $w.$ Let $\ell'_{vw}$ be the linking number of
$v$ and $w,$ excluding the weights around $v$ and $w.$

\begin{definition}[Semigroup Conditions]\label{sgcdef}The
\textit{semigroup condition} at $v$ in the direction of $e$ is
$$d_{ve}\in\N\langle \ell'_{vw}~|~w\text{ is a leaf in }\Delta_{ve}
\rangle.$$ We say that $\Delta$ \textit{satisfies the semigroup
conditions} if the semigroup condition for every node $v$ and
every adjacent edge $e$ is satisfied. Note that for an edge
leading to a leaf, the condition is trivially satisfied.
\end{definition}

To each leaf $w$ in $\Delta,$ associate a variable $Z_{w}.$  If $\Delta$
satisfies the semigroup conditions, then for each $v$ and $e$ as
above, there exist $\alpha_{vw}\in\N\cup\{0\}$ such that
$$d_{ve}=\sum_{w\text{ a leaf in }\Delta_{ve}}\alpha_{vw}\ell'_{vw}.$$ Then a monomial
$M_{ve}=\prod_{w}Z_{w}^{\alpha_{vw}},$ a product over leaves $w$ in $\Delta_{ve}$ with $\alpha_{vw}$ as
above,
is called an \textit{admissible monomial} for $e$ at $v.$  
If one associates the weight $\ell_{vw}$ to $Z_w,$ then for this
weight system, the so-called $v$-\textit{weighting}, $M_{ve}$ has weight $d_{v}=\prod_{e}d_{ve},$
where the product is taken over all edges $e$ adjacent to $v.$

\begin{definition}[Splice Diagram Equations]
Suppose $\Delta$ satisfies the semigroup conditions.  For each node $v$
and adjacent edge $e,$ choose an admissible monomial $M_{ve}.$  Let
$\delta_v$ denote the valency of the vertex $v.$  A a set of
\textit{splice diagram equations} for $\Delta$ is a set of equations of
the form
$$\left\{\sum_{e}a_{vie}M_{ve}=0~:~1\leq i\leq
\delta_v-2,~v\text{ a node in }\Delta\right\},$$ where for each $v,$
all maximal minors of the matrix $(a_{vie})$ have full rank. (One
can also add to each equation a convergent power series in the
$Z_w$ for which all of the terms have $v$-weight greater than $d_v.$  Since this
extension has no bearing upon the work herein, we omit it in
further discussion.)
\end{definition}

Each vertex $v\in\Gamma$ corresponds to an exceptional curve $E_v.$
Let $\mathbb{E}:=\bigoplus_{v\in\Gamma}\Z E_v.$  The intersection
pairing defines a natural injection
$\mathbb{E}\hookrightarrow\mathbb{E}^{*}=\text{Hom}(\mathbb{E},\Z),$
and the discriminant group is the finite abelian group
$D(\Gamma):=\mathbb{E}^{*}/\mathbb{E}.$ This group is isomorphic to $H_1(\Sigma, \Z).$ The order of $D(\Gamma)$ is
$\det(\Gamma):=\det(-C(\Gamma)),$ where
$C(\Gamma):\mathbb{E}\times\mathbb{E}\rightarrow\Z$ is the
intersection pairing.  There are induced symmetric pairings of
$\mathbb{E}\otimes\Q$ into $\Q$ and $D(\Gamma)$ into $\Q.$

Suppose $\Delta$ has $t$ leaves, and let $Z_1,\ldots,Z_t$ be the
associated variables. Neumann and Wahl define a faithful diagonal representation of $D(\Gamma)$ on $\C[Z_1,\ldots,Z_t].$ Let
$E_1,\ldots,E_t$ be the curves in $\Gamma$ corresponding to the
$t$ leaves of $\Delta,$ and let $e_j\in\mathbb{E}^{*}$ be the dual
basis element corresponding to $E_j.$  That is,
$e_j(E_k)=\delta_{jk}.$ Finally, for $r\in\Q,$ let $[r]$ denote
the image of the equivalence class of $r$ under the map
$\Q/\Z\hookrightarrow\C^{*}$ defined by $r\mapsto\exp(2\pi i r).$
Then the action of the discriminant group on the polynomial ring
$\C[Z_1,\ldots,Z_t]$ is generated by the action of the $e_j,~1\leq
j\leq t,$ which is defined by $e_j\cdot Z_k=[-e_j\cdot e_k]Z_k,$
$1\leq j,k\leq t.$

\begin{definition}[Congruence conditions]\label{ccdef}
Let $\Gamma$ be a graph for which the associated splice diagram $\Delta$
satisfies the semigroup conditions.  Then we say that $\Gamma$ satisfies
the \textit{congruence condition} at a node $v$ if one can choose an
admissible monomial for each adjacent edge $e$ such that all of
these monomials transform by the same character under the action of
$D(\Gamma).$  If this condition is satisfied for every node $v,$ then
$\Gamma$ \textit{satisfies the congruence conditions.}
\end{definition}  

We should mention here that Okuma gives a single condition that is
equivalent to the semigroup and congruence conditions together,
``Condition 3.3" of \cite{Okuma1}.  That this condition is
equivalent to the semigroup and congruence conditions is shown in
\cite{NW}. We will often say ``$\Gamma$ satisfies the semigroup and
congruence conditions", as opposed to ``$\Delta$ satisfies the
semigroup conditions and $\Gamma$ satisfies the congruence
conditions". Suppose a resolution graph $\Gamma$ satisfies the
semigroup and congruence conditions. Then, by a set of splice
diagram equations for $\Gamma,$ we mean equations as in Definition
\ref{sgcdef} such that for each $v,$ the admissible monomials
$M_{ve}$ transform equivariantly under $D(\Gamma).$  A resolution tree
$\Gamma$ is \textit{quasi-minimal} if any string in $\Gamma$ either
contains no $(-1)$-weighted vertex, or consists of a unique
$(-1)$-weighted vertex.

\begin{theorem}[\cite{NW}]\label{bignwthm} Suppose $\Gamma$ is quasi-minimal
and satisfies the semigroup and congruence conditions. Then a set
of splice diagram equations for $\Gamma$ defines an isolated complete
intersection singularity $\sing{Y},$ $D(\Gamma)$ acts freely on
$Y-\{0\},$ and the quotient $X:=Y/D(\Gamma)$ has an isolated normal
surface singularity and a resolution with dual resolution
graph $\Gamma.$ Moreover, $\sing{Y}\rightarrow \sing{X}$ is the
universal abelian cover.
\end{theorem}

We will use the next two propositions to check the congruence
conditions.

\begin{proposition}[\cite{NW}]\label{6.8} Let $\Gamma$ be a graph for
which the associated splice diagram $\Delta$ satisfies the
semigroup conditions.  Then the congruence conditions are
equivalent to the following: For every node $v$ and adjacent edge
$e$ in $\Delta,$ there is an admissible monomial
$M_{ve}=\prod_{w} Z_w^{\alpha_w}$ such that for
every leaf $w'$ in $\Delta_{ve},$
\begin{equation*}\left[\sum_{w\neq
w'}\alpha_w\frac{\ell_{ww'}}{\det(\Gamma)} -\alpha_{w'}e_{w'}\cdot
e_{w'}\right]=\left[\frac{\ell_{vw'}}{\det(\Gamma)}\right].\end{equation*}
\end{proposition}

\begin{remark}\label{nocheckleaves}It is easy to check, using the following proposition,
that this condition is always satisfied for an edge leading
directly to a leaf.\end{remark}
\begin{proposition}[\cite{NW}]\label{disc2}Suppose we have a string from a
leaf $w$ to an adjacent node $v$ in a resolution graph $\Gamma$ as in
the following diagram, with associated continued fraction $d/p.$ \begin{center}
\includegraphics[width=2in]{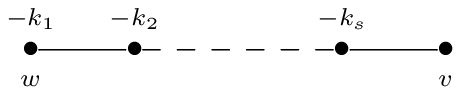}
\end{center} That is,
$$\frac{d}{p}=k_1-\dfrac{1}{k_2-\dfrac{1}{\ddots -\dfrac{1}{k_s}}}.$$ Then, if $d_v$ is the
product of weights at $v,$ $e_w\cdot
e_w=-d_v/(d^2\det(\Gamma))-p/d.$\end{proposition}

%%%%%%%%%%%%%%%%%%%%%%%%%%%%%%%%%%%%%%%%%%%%%%%%%%%%%%%%%%
%%%%%%%%%%%%%%%%%%%%% Section 3 %%%%%%%%%%%%%%%%%%%%%%%%%%%%%%%%
%%%%%%%%%%%%%%%%%%%%%%%%%%%%%%%%%%%%%%%%%%%%%%%%%%%%%%%%%%

\section{The resolution graph and splice diagram}\label{top}

Let $\{f(x,y)=0\}\subset\C^2$ define an analytically irreducible
plane curve with a singularity at the origin, and let
$X_{f,n}:=\{z^n=f(x,y)\}\subset\C^3.$ In \cite{MenNem}, Mendris
and N\'{e}methi prove that the link of $\xfn$ completely
determines the Newton/topological pairs of $f$ and the value of
$n,$ with two well-understood exceptions. In doing so, they give a
presentation of the construction of the resolution graph of $\xfn$
that is very useful
for our purposes. 
Section \ref{resgraph_sec} is a summary of the results we need
from Mendris and N\'{e}methi's work, and we use their notation
whenever possible. In section \ref{spldiag_sec}, we describe the associated
splice diagram.

It turns out that when $n=p_s=2,$ the resolution graph has a
structure that differs significantly from the general case. It is
referred to as the ``pathological case" or ``P-case" by Mendris
and N\'{e}methi, and we use this terminology as well. Some of the
computations must be done separately for the pathological case.

\subsection{Resolution graph}\label{resgraph_sec}

Suppose that $f$ has \textit{Newton pairs} $\{(p_k,q_k)~|~1\leq
k\leq s\}$ (see \cite{EN}, p. 49). They satisfy the following
properties: $q_1>p_1,$ $q_k\geq 1,~p_k\geq 2,\text{ and
}\gcd(p_k,q_k)=1$ for all $k.$ Define integers $a_k$ by $a_1=q_1,$
and
\begin{equation}\label{ak}a_{k}=q_{k}+a_{k-1}p_{k-1}p_{k},~2\leq k\leq s.\end{equation}
 The pairs
$\{(p_k,a_k)~|~1\leq k\leq s\},$ defined by Eisenbud and Neumann
in \cite{EN}, are referred to as the \textit{topological pairs} of
$f.$ These are the integers that appear in the splice diagram of
the link of the plane curve singularity defined by $f=0$ in $\C^2.$ Note
that $a_1>p_1,~a_k>a_{k-1}p_{k-1}p_k,~\text{and }\gcd(p_k,a_k)=1$
for all $k.$

The topological pairs $\{(p_k,a_k)~|~1\leq k\leq s\}$ are related to the \textit{Puiseux pairs} $\{(p_k,m_k)~|~1\leq k\leq s\}$ as follows: 
$a_1=m_1, \text{ and }a_k=m_k-m_{k-1}p_k+a_{k-1}p_{k-1}p_k,$ for $2\leq k\leq s. $  Furthermore, let $\bar{\beta_k},~0\leq k\leq s,$ be the generators of the semigroup associated to the plane curve singularity defined by $f$ (see \cite{Zariski}). Then we have $\bar{\beta}_0=p_1p_2\cdots p_s,~\bar{\beta}_k=a_k p_{k+1}\cdots p_s$ for $1\leq k\leq s-1,$ and $\bar{\beta}_s=a_s.$

\begin{figure}[h]
\centering
\includegraphics[width=2.5in]{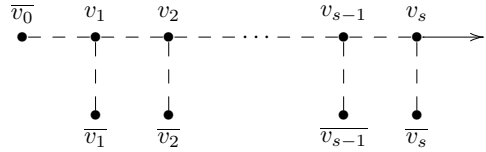}
\caption{Schematic form of
$\Gamma(\C^2,f),$ reproduced from \cite{MenNem}.}\label{schembres}
\end{figure}

By an embedded resolution of the germ of a function
$g:\sing{X}\rightarrow \sing{\C}$ we mean a resolution of the
singularity $\pi:\tilde{X}\rightarrow X$ such that
$\pi^{-1}(\{g=0\})$ is a divisor with only normal crossing singularities. We
also assume that no irreducible component of the exceptional set
$\pi^{-1}(0)$ intersects itself and that any two irreducible components have at most one intersection point. The minimal good embedded
resolution graph of $f:\sing{\C^2}\rightarrow\sing{\C}$ is a tree of rational
curves, denoted $\Gamma(\C^2,f).$ 
The construction of
the graph $\Gamma(\C^2,f)$ is well-known (e.g., \cite{BK}).
Reproducing the notation of Mendris and N\`{e}methi \cite{MenNem},
we consider this graph in a convenient schematic form (Figure
\ref{schembres}),
where the dashed lines represent strings of rational curves
(possibly empty) for which the self-intersection numbers are
determined by the continued fraction expansions of $p_k/q_k$ and
$q_k/p_k$ (see $\S$\ref{strings_sec} for details). 

There is an algorithm for constructing an embedded resolution
graph (not necessarily minimal) of the function $z:\xfn\rightarrow
(\C,0)$ from the graph $\Gamma(\C^2,f).$ 
Here, we follow the presentation in \cite{MenNem}, reproducing only what is necessary for
our purposes.  The output of this algorithm, without any
modifications by blow up or down, is referred to by Mendris and
N\'{e}methi as the \textit{canonical} embedded resolution graph of
$z$ in $\xfn,$ and is denoted $\Gamma^{can}(X_{f,n},z).$ The $n$-fold ``covering"
or ``graph projection" produced in the algorithm is denoted
$q:\Gamma^{can}(X_{f,n},z)\rightarrow\Gamma(\C^2,f).$
\begin{definition}[\cite{MenNem}]\label{allthestuff} Define positive integers $d_k,~h_k,~\h{k},~p'_{k},$ and $a'_{k}$ as
follows: $$\begin{array}{ll}\bullet\ d_k=(n,p_{k+1}p_{k+2}\cdots
p_s)& \text{for }0\leq k\leq s-1,\\
\bullet\ d_s=1;&\end{array}$$ and, for $1\leq k\leq s,$
$$\begin{array}{ll}\bullet\  h_k=(p_k,n/d_k),\quad &\bullet\ p'_k={p_k}/{h_k},\\
\bullet\ \h{k}=(a_k,n/d_k), &\bullet\  a'_k=a_k/\h{k}.
\end{array}$$\end{definition}%
If $w$ is a vertex in $\Gamma(\C^2,f),$ then all vertices in $q^{-1}(w)$
have the same multiplicity and genus, which we denote $m_w$ and
$g_w,$ respectively.

\begin{proposition}[\cite{MenNem}]\label{resgraph}
Let $q:\Gamma^{can}(X_{f,n},z)\rightarrow\Gamma(\C^2,f)$ be the ``graph projection"
mentioned above.  Then $\Gamma^{can}(X_{f,n},z)$ is a tree such that the following hold:\\

$\begin{array}{llll} (a) & \#q^{-1}(v_s)=1,&
\#q^{-1}(v_k)=h_{k+1}\cdots h_s, & (1\leq k\leq s-1)\\
& \#q^{-1}(\vbar{s})=\h{s},&
\#q^{-1}(\vbar{k})=\h{k}h_{k+1}\cdots h_s, & (1\leq k\leq s-1)\\
& \#q^{-1}(\vbar{0})=h_1\cdots h_s;& &
\end{array}$\\

$\begin{array}{lllll}
(b) & m_{v_k}&=&a'_k p'_k p'_{k+1}\cdots p'_s &(1\leq k\leq s),\\
& m_{\vbar{0}}&=& p'_1p'_2\cdots p'_s,&\\
& m_{\vbar{k}}&=&a'_k p'_{k+1}\cdots p'_s &(1\leq k\leq s-1),\\
& m_{\vbar{s}}&=&a'_s; &
\end{array}$\\

$\begin{array}{lllll}
(c) & g_{\vbar{k}}&=&0 &(0\leq k\leq s),\\
& g_{v_{k}}&=& (h_k-1)(\h{k}-1)/2 &(1\leq k\leq s).
\end{array}$\\

In particular, the link of $\xfn$ is a $\Q$HS if and only if
$(h_k-1)(\h{k}-1)=0$
for all $k,~1\leq k\leq
s.$
\end{proposition}

\begin{figure}[h]\begin{center}
\leavevmode
\includegraphics[width=\textwidth]{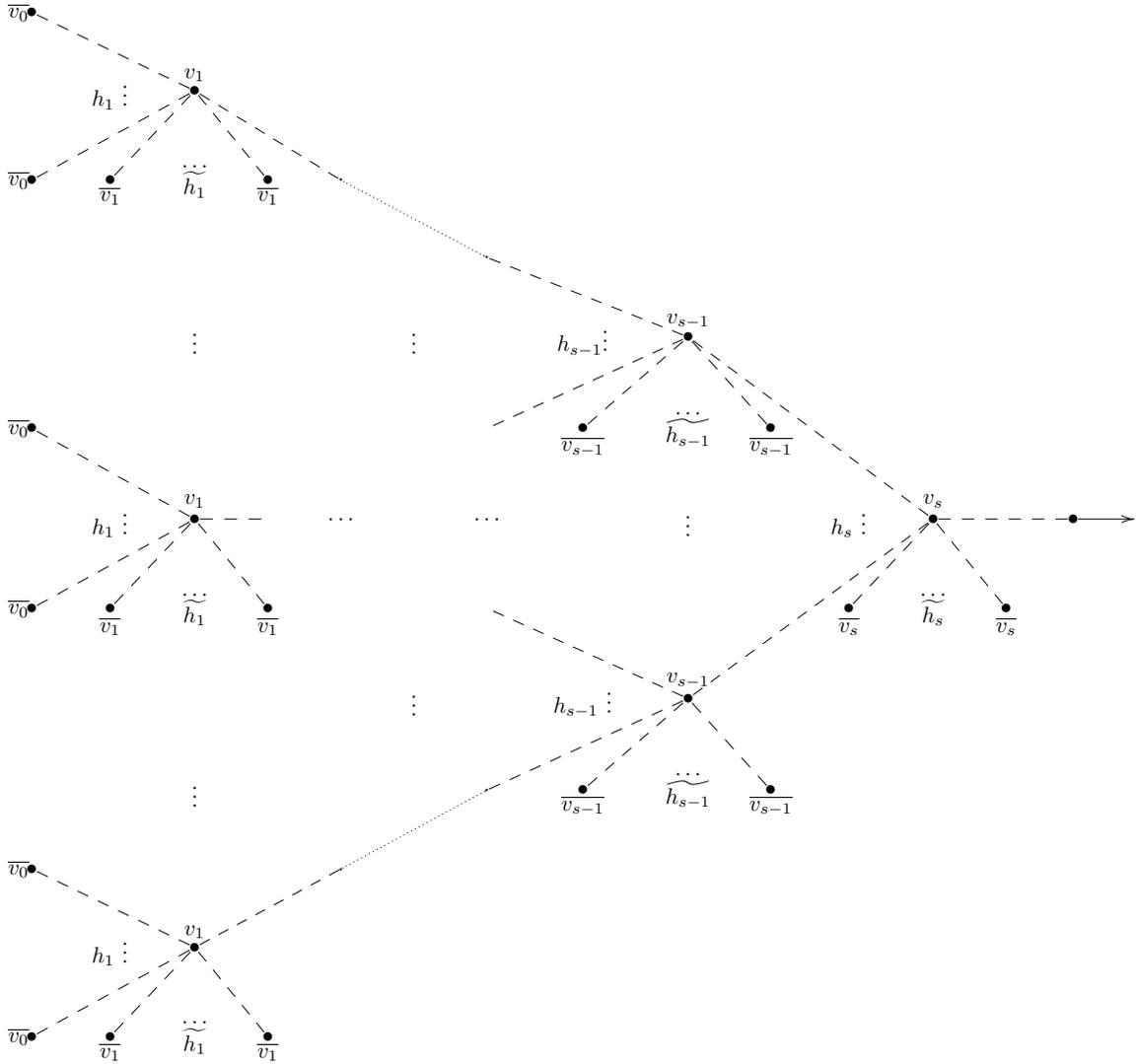}
\caption{Schematic form of $\Gamma^{can}(X_{f,n},z),$ reproduced from
\cite{MenNem}.}\label{biggraph}\end{center}
\end{figure}

The schematic form of $\Gamma^{can}(X_{f,n},z)$ is displayed in Figure
\ref{biggraph}, which is reproduced from \cite{MenNem}.
Abusing notation, we have labelled any vertex in $q^{-1}(v_k)$
(respectively, $q^{-1}(\vbar{k})$) with $v_k$ (respectively,
$\vbar{k}$). The dashed lines represent strings of vertices that
are not necessarily minimal. By the construction, each string must
contain at least as many vertices as its image in $\Gamma(\C^2,f).$ A vertex
is called a \textit{rupture vertex} if either it has positive
genus or it is a node. Note that any rupture vertex of $\Gamma^{can}(X_{f,n},z)$
must be in $q^{-1}(v_k)$ for some $k.$

\subsubsection*{Certain subgraphs of $\Gamma^{can}(X_{f,n},z)$ and their determinants}
Let $w$ be a vertex in $\Gamma(\C^2,f),$ and let $v'$ be any vertex in
$q^{-1}(w).$ If $w=v_k$ for some $k,~1\leq k\leq s-1,$ then the
shortest path from $v'$ to the arrowhead of $\Gamma^{can}(X_{f,n},z)$ contains at
least one rupture vertex, and the rupture vertex along that path
which is closest to $v'$ is a vertex $v''\in q^{-1}(v_{k+1}).$
Define $\Gamma(v')$ to be the subgraph of $\Gamma^{can}(X_{f,n},z)$ consisting of the
string of vertices between $v'$ and $v'',$ not including $v'$ and
$v''.$ If $w=v_s,$ then the shortest path from $v'$ to the
arrowhead is a string; let $\Gamma(v')$ be this string, not including
$v'.$ Finally, if $w=\vbar{k},~0\leq k\leq s,$ let $v''$ be the
rupture vertex that is closest to $v'$ on the shortest path from
$v'$ to the
arrowhead. 
Define $\Gamma(v')$ to be the subgraph of consisting of the string of
vertices from $v'$ to $v'',$ including $v'$ but not $v''.$ Up to
isomorphism, none of these strings depend upon the choice of $v'$
in $q^{-1}(w),$ so whenever the particular vertex $v'$ does not
matter, we will simply denote them $\Gamma(w).$

Fix an integer $k,~1\leq k\leq s,$ and fix a vertex $v'$ in
$q^{-1}(v_k).$ Consider the collection of connected subgraphs that
make up $\Gamma^{can}(X_{f,n},z)-\{v'\}.$ There are $\h{k}$ isomorphic components
that are strings of isomorphism type $\Gamma(\vbar{k}).$ There is one
connected subgraph that contains the arrowhead; denote this
subgraph $\Gamma_{A}(v').$ The $h_k$ remaining components are all
isomorphic. Let $\Gamma_{-}(v')$ denote any of these isomorphic
subgraphs. Again, whenever the particular choice of $v'$ is
unimportant, we use $\gm{k}$ instead of $\Gamma_{-}(v'),$ and
$\ga{k}$ instead of $\Gamma_{A}(v').$ Note that
$\gm{1}=\Gamma(\vbar{0})$ and $\ga{s}=\Gamma(v_s).$ We should also point
out that the subgraphs $\ga{k}$ do not appear in \cite{MenNem}; in
particular, $\ga{k}$ is not the same as their $\Gamma_{+}(v_k).$

For any resolution graph $\Gamma,$ let $\det(\Gamma):=\det(-C),$ where $C$
is the intersection matrix of the exceptional curves in $\Gamma.$ If
$\Gamma$ is empty, then we define $\det(\Gamma)$ to be $1.$ Nearly all of
the determinants of the subgraphs defined above are explicitly
computed by Mendris and N\'{e}methi in \cite{MenNem}, and those
that are not can be computed by the same method. 
\begin{lemma}[\cite{MenNem}]\label{detmaxstrings} For any
$w$ in $\Gamma(\C^2,f)$ as above, let $D(w):=\det(\Gamma(w)).$ Then
$$\begin{array}{lll}
D(\vbar{0})&=&a'_1,\\
D(\vbar{k})&=&p'_k,\text{ for } 1\leq k\leq s,\\
D(v_s)&=&n/(h_s\h{s}),\\
D(v_k)&=&n q_{k+1}/(d_{k-1}\h{k}\h{k+1}),\text{ for }1\leq k\leq
s-1.
\end{array}$$
\end{lemma}
It follows from the construction of $\Gamma^{can}(X_{f,n},z)$ that if $D(v_s)=1,$
this indicates that $\Gamma(v_s)$ is empty, and the arrowhead in
$\Gamma^{can}(X_{f,n},z)$ is connected directly to the unique vertex in $q^{-1}(v_s).$ 
\begin{lemma}[\cite{MenNem}]\label{dm} Let $\dm{k}:=\det(\gm{k}),~1\leq k\leq s.$ If $s\geq
2,$ then for $2\leq k\leq s,$
$$\frac{\dm{k}}{a'_k}=(a'_{k-1})^{h_{k-1}-1}(p'_{k-1})^{\h{k-1}-1}
\left[\frac{\dm{k-1}}{a'_{k-1}}\right]^{h_{k-1}}.
$$
\end{lemma} 

The method used to prove Lemma \ref{dm} can be suitably modified
to prove the next two lemmas. The computation is straightforward,
so we omit the proof.
\begin{lemma}\label{dalemma} Assume $s\geq 2,$ and let $\da{k}:=\det(\ga{k}),~1\leq k\leq s.$
Let $A_k$ be defined recursively by
$A_{s-1}=a_{s-1}p_{s-1}p'_s+q_s,$ and, for $1\leq k\leq s-2,$
$$A_k=a_{k}p_{k}p'_{k+1} A_{k+1} +q_{k+1}a_{k+2}\cdots a_s.$$
Then
\begin{equation*}\da{k}=\frac{n A_k \left\{\prod_{j=k+1}^{s}
(p'_j)^{\h{j}-1}\dm{j}^{h_j-1}\right\}}{\HH{k}d_ka_{k+1}\cdots
a_s},\text{ for }1\leq k\leq s-1.
\end{equation*}
\end{lemma}

\begin{lemma}\label{det} The determinant of $\Gamma^{can}(X_{f,n},z)$ is given by
\begin{equation*}\det(\Gamma^{can}(X_{f,n},z))=(a'_s)^{h_s-1}(p'_s)^{\h{s}-1}\left[\frac{\dm{s}}{a'_s}\right]^{h_s}.
\end{equation*}
\end{lemma}

A minimal good embedded resolution graph of $z$ in $\xfn,$ denoted
$\Gamma^{min}(X_{f,n},z),$ is obtained from $\Gamma^{can}(X_{f,n},z)$ by repeatedly blowing down any
rational $(-1)$-curves for which the
corresponding vertex has valency one or two. By dropping the
arrowhead and multiplicities of $\Gamma^{min}(X_{f,n},z)$ and then blowing down any
appropriate rational $(-1)$-curves, we obtain a
minimal good resolution graph of $\xfn,$ denoted $\Gamma^{min}(X_{f,n}).$

\begin{proposition}[\cite{MenNem}]\label{nodesurvival}All of the rupture vertices
in $\Gamma^{can}(X_{f,n},z)$ survive as rupture vertices in $\Gamma^{min}(X_{f,n},z).$  That is, they
are not blown down in the minimalization process, and after
minimalization, they are still rupture vertices.
\end{proposition}

\begin{proposition}[\cite{MenNem}]\label{pathcase}Assume that by deleting the arrowhead of $\Gamma^{min}(X_{f,n},z)$
we obtain a non-minimal graph.  This situation can happen if and
only if $n=p_s=2.$ In this case, the link is a $\Q$HS and $\Gamma^{min}(X_{f,n},z)$
has the following schematic form, with $e\geq 3.$%
\begin{center}
\includegraphics[width=3in]{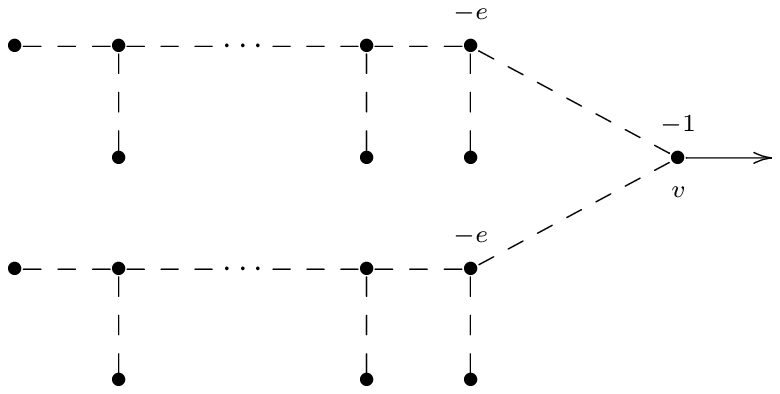}
\end{center}%
The minimal resolution graph $\Gamma^{min}(X_{f,n})$ is obtained from $\Gamma^{min}(X_{f,n},z)$ by deleting
the arrowhead and blowing down $v.$
\end{proposition}
Propositions \ref{nodesurvival} and \ref{pathcase} imply that all
of the nodes in $\Gamma^{can}(X_{f,n},z)$ remain nodes in the minimal good
resolution graph of $\xfn$ except in the case $n=p_s=2.$ We refer
to $n=p_s=2$ as the \textit{pathological case}, and it is treated
separately in what follows.

%%%%%%%%%%%%%%%%%%%%%%%%%%%%%%%%%%%%%%%%%%%%%%%%%%%%%%%%%%%%%%%%%
%%%%%%%%%%%%%%%%%%%%%%%%%%%%%%%%%%%%%%%%%%%%%%%%%%%%%%%%%%%%%%%%%

\subsection{Splice diagram}\label{spldiag_sec}

From now on, we assume that the link of $\xfn$ is a $\qhs.$ That is, for each $k,~1\leq k\leq s,$ either $h_k$ or $\h{k}$ is equal to $1.$ 
One complication that arises is that certain strings in $\Gamma^{can}(X_{f,n},z)$
may completely collapse upon minimalization. Therefore, if we use
the minimal good resolution graph $\Gamma^{min}(X_{f,n})$ in what follows, we
would constantly need to note that certain strings may be empty,
and more importantly, that certain leaves in the splice diagram
may not be present. We will avoid this by using the splice diagram associated to $\Gamma^{can}(X_{f,n}),$ the graph that results from deleting the arrowhead and multiplicities in $\Gamma^{can}(X_{f,n},z).$
We could easily use a quasi-minimal modification of $\Gamma^{can}(X_{f,n}),$ and the computation of the splice diagram would 
not change.  Therefore, we can apply Theorem \ref{bignwthm} to $\Gamma^{can}(X_{f,n}).$

\subsubsection*{Splice diagram in the general case}
Assume we are not in the pathological case, and let $\Delta_{f,n}$ be
the splice diagram associated to $\Gamma_{f,n}:=\Gamma^{can}(X_{f,n}).$ 
If a vertex $v$ in $\Gamma_{f,n}$ is in
$q^{-1}(v_k)$ (respectively, $q^{-1}(\vbar{k})$), we say that $v$ is ``of type $v_k$" (respectively,  $\vbar{k}$). 
We use the same terminology for the
corresponding vertices of $\Delta_{f,n}.$

Consider a node $v$ of type $v_k,~1\leq k\leq s,$ in $\Gamma_{f,n}.$
In general, there are $h_k+\h{k}+1$ edges adjacent to $v$: $\h{k}$
edges that lead to strings of (isomorphism) type $\Gamma(\vbar{k}),$
$h_k$ edges that lead to subgraphs of type $\gm{k},$ and $1$ edge
that leads towards a subgraph of type $\ga{k}.$ 
The
corresponding pieces of $\Delta_{f,n}$ associated to the subgraphs of
type $\gm{k}$ and $\ga{k}$ are denoted $\delm{k}$ and $\dela{k},$
respectively. Recall that $\gm{1}=\Gamma(\vbar{0}),$ and
$\ga{s}=\Gamma(v_s),$ and keep in mind that $\Gamma(v_s)$ may be empty.

\begin{figure}[h]
\centering
\includegraphics[width=2in]{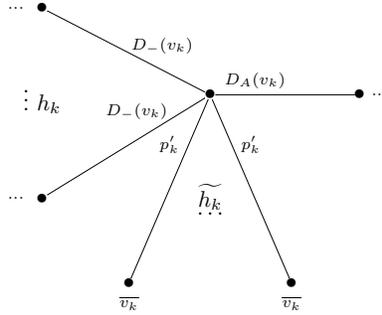}
\caption{Splice diagram at a node of type $v_k,~2\leq k\leq
s-1.$}\label{sdk}
\end{figure}

The weights of the splice
diagram $\Delta_{f,n}$ are given by Lemmas \ref{detmaxstrings},
\ref{dm}, and \ref{dalemma}. At a node of type $v_k$ in
$\Delta_{f,n},$ the weights on the $\h{k}$ edges that lead to leaves of
type $\vbar{k}$ are $D(\vbar{k})=p'_k;$ the weights on the $h_k$
edges connected to subgraphs of type $\delm{k}$ are $\dm{k};$ and
the weight on the single edge connected to the subgraph of type
$\dela{k}$ is $\da{k}$ (see Figure \ref{sdk}).

\subsubsection*{The pathological case} 
\begin{figure}[h]\centering
\includegraphics[width=4in]{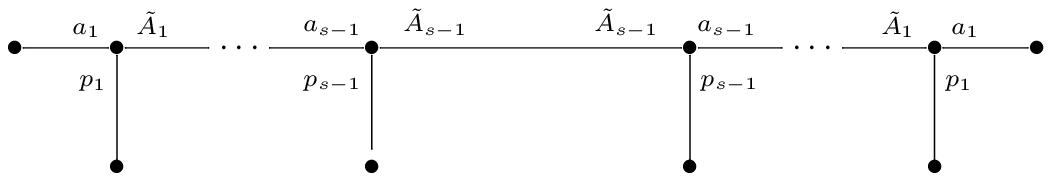}
\caption{Splice diagram for the pathological case.}\label{pathsd}
\end{figure}

For this case ($n=p_s=2$), it is more convenient to use the splice diagram 
associated to the minimal resolution graph $\Gamma^{min}(X_{f,n})$ (see Figure \ref{pathsd}).
Here, $h_s=2,$ hence $n/h_s=n/d_{s-1}=1.$ Then, by
definition $h_k=\h{k}=1$ for $1\leq k\leq s-1,$ and $\h{s}=1$
since $\gcd(p_s,a_s)=1.$ The link is a $\qhs,$ and the only string
of type $\Gamma(\vbar{k})$ that collapses completely in $\Gamma^{min}(X_{f,n},z)$ is
$\Gamma(\vbar{s})$ (Proposition \ref{pathcase}). The graph $\Gamma^{min}(X_{f,n})$ has a total of $2(s-1)$ nodes: two of
type $v_k$ for each $k,~1\leq k\leq s-1.$ Each of these nodes has
valency three.

Since the determinant of a resolution tree remains constant
throughout the minimalization process, the weights of the splice diagram associated to $\Gamma^{min}(X_{f,n})$ can be
determined from Lemmas \ref{detmaxstrings}, \ref{dm}, and
\ref{dalemma}. Since $\HH{k}=1$ for $1\leq k\leq s-1,$ we have
$\dm{k}=a_k$ for $2\leq k\leq s.$ Define integers $\tilde{A_k}$ as
follows:
\begin{eqnarray*}\tilde{A}_{k}&:=&a_s-a_{k}p_{k}p_{k+1}^2\cdots
p_{s-1}^2,\text{ for }1\leq k\leq
s-2,\text{ and}\\
\tilde{A}_{s-1}&:=&a_s-a_{s-1}p_{s-1}.\end{eqnarray*}  It is easy
to check that $\da{k}=\tilde{A}_k$ for $1\leq k\leq s-1.$

%%%%%%%%%%%%%%%%%%%%%%%%%%%%%%%%%%%%%%%%%%%%%%%%%%%%%%%%%%
%%%%%%%%%%%%%%%%%%%%% Section 4 %%%%%%%%%%%%%%%%%%%%%%%%%%%%%%%%
%%%%%%%%%%%%%%%%%%%%%%%%%%%%%%%%%%%%%%%%%%%%%%%%%%%%%%%%%%

\section{The semigroup conditions}\label{sgc_sec}

In this section, we discuss the semigroup conditions for the
splice diagram $\Delta_{f,n}.$  Throughout this section, we assume
that we are not in the pathological case. For a node $v$ of type
$v_k$ in $\Delta_{f,n},$ $1\leq k\leq s,$ there are at most two
inequivalent semigroup conditions to check: one for an edge that
leads to a subdiagram of type $\delm{k}$ (nontrivial for $2\leq
k\leq s$), and one for the edge that leads to a subdiagram of type
$\dela{k}$ (nontrivial for $1\leq k\leq s-1$). Clearly, for a
fixed $k,$ the semigroup conditions are equivalent for any node
$v$ of type $v_k.$

\subsection*{Semigroup conditions in the direction of $\delm{k}$ }
\begin{lemma}\label{genclaim} Let $v$ be a node of type $v_k, ~2\leq k\leq s,$
and let $w_j$ be a leaf of type $\vbar{j}$ in $\Delta_{-}(v),~0\leq j\leq
k-1.$ Then
$$\ell'_{vw_j}=\left\{\begin{array}{ll}
(\dm{k}/a'_k)p'_1\cdots p'_{k-1}&\text{ for }j=0\\
(\dm{k}/a'_k)a'_jp'_{j+1}\cdots p'_{k-1}&\text{ for }1\leq j\leq k-2\\
(\dm{k}/a'_k)a'_{k-1} & \text{ for }j=k-1.
\end{array}\right.$$\end{lemma}

\begin{proof} We prove this by induction on $k.$
For $k=2,$ the lemma is true,
since if $v$ is a node of type $v_2,$
\begin{eqnarray*}\ell'_{vw_0}&=&(a'_1)^{h_1-1}(p'_1)^{\h{1}},\\
\ell'_{vw_1}&=&(a'_1)^{h_1}(p'_1)^{\h{1}-1},\text{ and }\\
\dm{2}/a'_2&=&(a'_1)^{h_1-1}(p'_1)^{\h{1}-1}.
\end{eqnarray*} 

Now assume the lemma is true for $k=i-1;$ we show that it is true
for $k=i.$ 
Fix a node $v$ of type $v_i,$ and (abusing notation), let
$v_{i-1}$ denote the unique node of type $v_{i-1}$ in $\Delta_{-}(v).$
For $0\leq j\leq i-2,$ any leaf of type $\vbar{j}$ in $\Delta_{-}(v)$ is
in one of the subdiagrams of type $\delm{i-1}.$ Thus (refer to
Figure \ref{fig3_1})
$$\ell'_{vw_j}=\left\{\begin{array}{ll}
\dm{i-1}^{h_{i-1}-1}(p'_{i-1})^{\h{i-1}}\ell'_{v_{i-1}w_j}&\text{
for }0\leq j\leq
i-2,\\
\dm{i-1}^{h_{i-1}}(p'_{i-1})^{\h{i-1}-1}&\text{ for }j=i-1.
\end{array}\right.$$
\begin{figure}[h]
\centering
\includegraphics[width=4in]{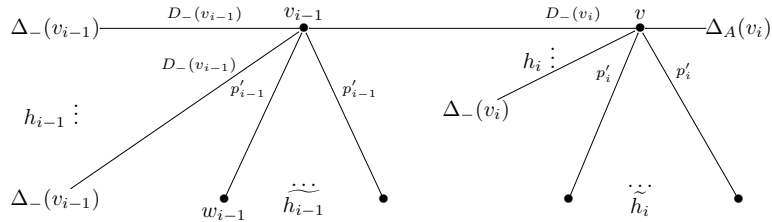}
 \caption{Relevant portion of $\Delta_{f,n}$ at a node $v$ of type
$v_i.$} \label{fig3_1}
\end{figure}

By Lemma \ref{dm}, we have
$\frac{\dm{i}}{a'_i}=(p'_{i-1})^{\h{i-1}-1}\dm{i-1}^{h_{i-1}-1}\cdot\frac{\dm{i-1}}{a'_{i-1}}.$
Applying this fact and the induction hypothesis yields the desired result.\end{proof}

\begin{proposition}\label{sgcleft} At a node of type $v_k,~2\leq k\leq s,$ the
semigroup condition in the direction of any of the $h_k$ edges
that lead to a subdiagram of type $\delm{k}$ is equivalent to
\begin{equation}\label{sgcminus}a'_k\in\N\langle a'_{k-1},~p'_1p'_2\cdots
p'_{k-1},~a'_{j}p'_{j+1}\cdots p'_{k-1},~1\leq j\leq
k-2\rangle.\end{equation} Furthermore, if $\h{k}=1,$ this
condition is automatically satisfied.
\end{proposition}

\begin{proof} Fix a node $v$ of type $v_k$ in $\Delta_{f,n}.$ By
Definition \ref{sgcdef}, the condition is
$$\dm{k}\in\N\langle \ell'_{vw}~|~w\text{ is a leaf in }\Delta_{-}(v)
\rangle.$$  
The leaves in $\Delta_{-}(v)$ are of type $\vbar{j},$ for $j$ such that
$0\leq j\leq k-1.$ Hence, there are $k$ generators for the
semigroup in question, namely, $\ell'_{vw_j},~0\leq j\leq k-1,$
where $w_j$ denotes any leaf in $\Delta_{-}(v)$ of type $\vbar{j}.$   The
first statement of the Proposition follows from Lemmas \ref{dm}
and \ref{genclaim}, since $\dm{k}$ and all generators of the
semigroup are divisible by $\dm{k}/a'_k$.

The second statement follows from \cite{zhs}, Proposition 8.1.
\end{proof}

\subsection*{Semigroup conditions in the direction of $\dela{k}$ }

\begin{figure}[h]
\centering
\includegraphics[width=\textwidth]{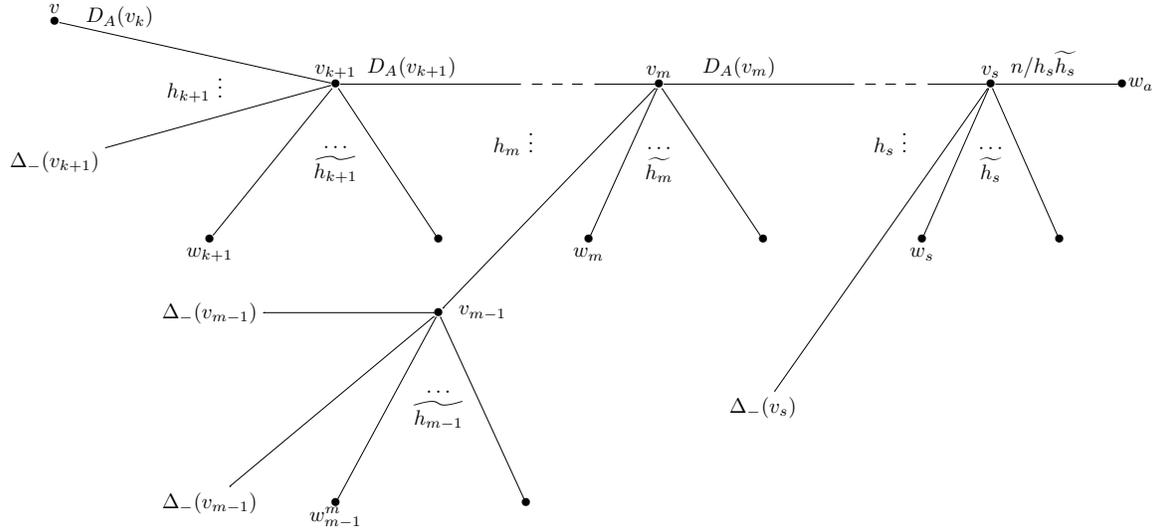}
\caption{Relevant portion of $\Delta_{f,n}$ at a node $v$ of type $v_k.$}
\label{fig3_2}
\end{figure}

Fix an integer $k,~1\leq k\leq s-1,$ and fix a node $v$ of type
$v_k.$ 
By definition, the semigroup condition is $\da{k}\in\mathrm{R}_k,$
where
$$\mathrm{R}_k:=\N\langle \ell'_{vw}~|~w\text{ is a leaf in }\Delta_{A}(v)
\rangle.$$ Refer to Figure \ref{fig3_2} for what follows.
There is at least one leaf $w_s$ in $\Delta_{A}(v)$ of type $\vbar{s}$
connected to $v_s$ (the unique node of type $v_s$), and if
$n/h_s\h{s}\neq 1,$ there is a leaf $w_a$ resulting from the
string $\Gamma(v_s)$ in $\Gamma_{f,n}.$ These contribute $\ell'_{vw_s}$ and
$\ell'_{vw_a}$ as generators of $\mathrm{R}_k.$

Next, travel along the shortest path from $v$ to $v_s.$  If
$k<s-1,$ this path contains one node of type $v_m,$ for each $m$
such that $k+1\leq m\leq s-1.$ Since there can be no confusion
here, we will simply refer to the nodes along this path as $v_m.$
Each of these nodes is directly connected to at least one leaf
$w_m$ of type $\vbar{m}.$ Each such leaf contributes the generator
$\ell'_{vw_m}$ to $\mathrm{R}_k.$ If $h_i=1$ for $k+1 \leq i\leq s,$
there are no other types of leaves in $\Delta_{A}(v),$ and we have listed
all the generators of $\mathrm{R}_k.$

For each $m$ such that $h_m\neq 1,~k+1\leq m\leq s,$ there are
more generators for $\mathrm{R}_k,$ namely $\ell'_{vw}$ for each type
of leaf $w$ in $\delm{m}.$  There are $m$ such different types of
leaves: type $\vbar{j},$ for $j$ such that $0\leq j\leq m-1.$ Let
$w_j^m$ be a leaf of type $\vbar{j}$ in $\delm{m}.$ 
Then the generators of the semigroup $\mathrm{R}_k$ are:
$$\left\{ \begin{array}{ll}
\ell'_{vw_m}, & k+1\leq m\leq s,\\
\ell'_{vw_j^m}, & 0\leq j\leq m-1, \text{for all }m\text{ such
that
}k+1\leq m\leq s \text{ and }h_m\neq 1\\
\ell'_{vw_a} & (\text{absent if }n/\HH{s}= 1)
\end{array}\right\}.$$

\begin{proposition}\label{steps1and2} Suppose $h_s>1.$ Then the semigroup conditions imply that
$h_s=p_s$ and $\HH{s-1}=1.$ \end{proposition}

\begin{proof} Note that since the link is a $\qhs,$
$h_s>1$ implies $\h{s}=1.$ Let $v$ be a node of type $v_{s-1},$
and consider the semigroup condition at $v$ in the direction of
$\Delta_A(v)$: $\da{s-1}$ is in the semigroup $\mathrm{R}_{s-1}.$ The
generators of $\mathrm{R}_{s-1}$ are $\ell'_{vw_s},$
$\ell'_{vw_j^s},~0\leq j\leq s-1,$ and $\ell'_{vw_a}$ (absent if
$n/h_s=1$).  
\begin{figure}[h] \centering
\includegraphics[width=4in]{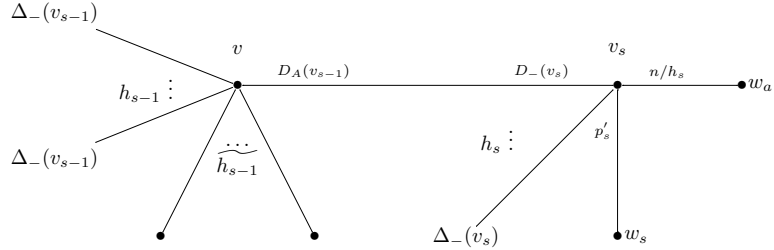}
\caption{Splice diagram $\Delta_{f,n}$ for $\h{s}=1.$}\label{fig_j}
\end{figure}
It is easy to check (see
Figure \ref{fig_j}) that
\begin{eqnarray*}
\ell'_{vw_s}&=&(n/h_s)\dm{s}^{h_s-1},\\
\ell'_{vw_a}&=&p'_s\dm{s}^{h_s-1},\text{ and }\\
\ell'_{vw_j^s}&=&(n/h_s)p'_s\dm{s}^{h_s-2}\ell'_{v_sw_j^s}.
\end{eqnarray*} By Lemma \ref{dalemma}, since $d_{s-1}=h_s$ and
$a'_s=a_s,$
$$\da{s-1}=\frac{nA_{s-1}\dm{s}^{h_s-1}}{\HH{s-1}h_sa_s},$$ where
$A_{s-1}=a_{s-1}p_{s-1}p'_s+q_s=a_s-a_{s-1}p_{s-1}(p_s-p'_s).$
Note that $n/(\HH{s-1}h_s)$ and $\dm{s}^{h_s-1}/a_s$ are both
integers in this case. Since $p_s>p'_s=p_s/h_s,$
$$\frac{n}{\HH{s-1}h_s}[a_s-a_{s-1}p_{s-1}(p_s-p'_s)]
<\frac{n}{\HH{s-1}h_s}a_s\leq \frac{n}{h_s}a_s,$$ and therefore
$\da{s-1}<\ell'_{vw_s}.$ Hence we can forget about the generator
$\ell'_{vw_s},$ since it is too large.

By Lemma \ref{genclaim},
$$\ell'_{v_sw_j^s}=\left\{\begin{array}{ll}
p'_1\cdots p'_{s-1}\cdot \dm{s}/a_s&\text{ for }j=0\\
a'_jp'_{j+1}\cdots p'_{s-1}\cdot \dm{s}/a_s&\text{ for }1\leq j\leq s-2\\
a'_{s-1}\cdot \dm{s}/a_s & \text{ for }j=s-1. \end{array}\right.$$
So, all generators of $\mathrm{R}_{s-1}$ and $\da{s-1}$ are divisible
by $\dm{s}^{h_s-1}/a_s,$ and the semigroup condition is equivalent
to the following: $n/(\HH{s-1}h_s)A_{s-1}$ is in the semigroup
generated by
\begin{equation}\label{newgens}\left\{\frac{n}{h_s}p'_1\cdots
p'_s,~\frac{n}{h_s}a'_jp'_{j+1}\cdots p'_s~:~1\leq j\leq
s-1,~a_sp'_s~(\text{Absent if }
\frac{n}{h_s}=1)\right\}.\end{equation} All of the generators of
this semigroup are divisible by $p'_s.$ Therefore, the semigroup
condition implies that $p'_s$ divides
$n/(\HH{s-1}h_s)[a_s-a_{s-1}p_{s-1}(p_s-p'_s)].$ Suppose $p'_s>1.$
Since $p'_s$ divides $p_s-p'_s,$ and $(a_s,p_s)=1,$ this implies
that $p'_s$ divides $n/(\HH{s-1}h_s).$ This is impossible, since
by definition $p'_s=p_s/(n,p_s),$ and thus $(p'_s,n)=1.$ Therefore
we must have $p'_s=1.$ Since $p'_s=p_s/h_s,$ we have shown that
the semigroup conditions imply $h_s=p_s.$

Now we show that the semigroup conditions imply $\HH{s-1}=1.$ Note
that if $n/h_s=1,$ this is automatically true by definition of
$h_i$ and $\h{i}.$ Therefore, assume that $n/h_s\neq 1.$ Observe
that all of the generators in (\ref{newgens}) are divisible by
$n/h_s$ except for $a_s.$ Therefore, if the semigroup condition is
satisfied, there exist $M$ and $N$ in $\N\cup\{0\}$ such that
$$n/(\HH{s-1}h_s)[a_s-a_{s-1}p_{s-1}(p_s-1)]=Ma_s+N n/h_s.$$
Hence,
\begin{eqnarray*}\left(n/(\HH{s-1}h_s)-M\right)a_s&=&N n/h_s
+n/(\HH{s-1}h_s)a_{s-1}p_{s-1}(p_s-1)\\
&=&n/h_s(N+a'_{s-1}p'_{s-1}(p_s-1)).
\end{eqnarray*} Since $(n,a_s)=1$ by assumption, this implies that
$n/h_s\neq 1$ divides $\frac{n}{\HH{s-1}h_s}-M.$ But we have
$$0<\frac{n}{\HH{s-1}h_s}-M\leq \frac{n}{\HH{s-1}h_s}\leq\frac{n}{h_s}.$$
Therefore, the only possibility is $n/h_s=n/(\HH{s-1}h_s)-M,$
i.e., $M=0$ and $\HH{s-1}=1.$\end{proof}

\begin{lemma}\label{induction} Assume $s\geq 3,$ and that
$\HH{s-1}=1.$  Then the semigroup conditions imply that $\HH{k}=1$
for $1\leq k\leq s-2.$
\end{lemma}

\begin{proof}  We prove this by strong downward induction on $k.$
First we show that the semigroup conditions imply that
$\HH{s-2}=1.$ By Proposition \ref{resgraph}(a), there are $h_s$
nodes of type $v_{s-2};$ let $v$ be any such node. We will show
that the semigroup condition for $v$ in the direction of $\Delta_{A}(v)$
cannot be satisfied if $\HH{s-2}\neq 1.$ 

Let $\tilde{A}_i=a_s-a_ip_ip_{i+1}^2\cdots
p_{s-1}^2(p_s-p'_s),~1\leq i\leq s-2.$ By Lemma \ref{dalemma},
$$\da{s-1}=\left\{\begin{array}{ll}
n(p_s)^{\h{s}-1} & \text{for } h_s=1\\
\frac{nA_{s-1}\dm{s}^{h_s-1}}{h_sa_s}& \text{for } h_s>1,
\end{array}\right.$$ and
$$
\da{s-2}=\left\{\begin{array}{ll}
\frac{n}{\HH{s-2}}(p_s)^{\h{s}-1}& \text{for } h_s=1\\
\frac{n\tilde{A}_{s-2}\dm{s}^{h_s-1}}{\HH{s-2}h_sa_s}& \text{for }
h_s>1.\end{array}\right.
$$
The generators of $\mathrm{R}_{s-2}$ are
$$\left\{\begin{array}{lcll} \ell'_{vw_{s-1}}&=&\da{s-1},&\\
\ell'_{vw_s}&=& n/(\HH{s})p_{s-1}\dm{s}^{h_s-1}
{p'_s}^{\h{s}-1},&\\
\ell'_{vw_j^s}&=& n/(\HH{s})p_{s-1}\dm{s}^{h_s-2}
{p'_s}^{\h{s}}\ell'_{v_sw_j^s},&0\leq j\leq
s-1,\\
\ell'_{vw_a}&=&p_{s-1}\dm{s}^{h_s-1} {p'_s}^{\h{s}}&.
\end{array}\right\},$$ although the $\{\ell'_{vw_j^s}\}_{j=0}^{s-1}$ are absent if
$h_s=1,$ and $\ell'_{vw_a}$ is absent if $n/\HH{s}=1.$

We will consider two separate cases: (i) $h_s=1,$ and (ii)
$h_s>1.$

Case (i). If $h_s=1,$ it is easy to see that if $\HH{s-2}\neq 1,$
then $\ell'_{vw_{s-1}}>\da{s-2}.$ Then, since $\da{s-2}$ and every
generator of the semigroup are divisible by $(p_s)^{\h{s}-1},$ the
semigroup condition is equivalent to: $n/(\HH{s-2})$ is in the
semigroup generated by $p_{s-1}n/\h{s}$ and
$p_{s-1}p_s~(\text{absent if }n/\h{s}=1).$  Thus the semigroup
condition implies that $n/(\HH{s-2})$ is divisible by $p_{s-1},$
which is impossible since $h_{s-1}=(n,p_{s-1})=1.$ Therefore, we
must have $\HH{s-2}=1.$ (Note that the argument is valid even if
$n/\h{s}=1$ or $\h{s}=1.$)

Case (ii). For $h_s>1,$ the proof that $\HH{s-2}$ must be $1$ is
nearly identical to the proof of Proposition \ref{steps1and2}, so
we just give the outline here. Recall that the semigroup
conditions imply that $p'_s=1$ in this case. Furthermore, we can
assume that $n/h_s\neq 1,$ since otherwise the Lemma is trivially
true by definition of $h_i$ and $\h{i}.$

Dividing $\da{s-2}$ and all the generators of $\mathrm{R}_{s-2}$
by $\dm{s}^{h_s-1}/a_s,$
we see that the semigroup condition for $v$ in the direction of
$\Delta_{A}(v)$ implies that $n/(\HH{s-2}h_s)\tilde{A}_{s-2}$ is in the
semigroup generated by $a_sp_{s-1}$ and a collection of positive
integers that are divisible by $n/h_s.$
The semigroup condition implies that there exist $M$ and $N$ in
$\N\cup\{0\}$ such that
$$n/(\HH{s-2}h_s)\tilde{A}_{s-2}=Ma_s
p_{s-1}+Nn/h_s.$$ Just as in the proof of Proposition
\ref{steps1and2}, we see that we must have $M=0$ and $\HH{s-2}=1.$
Thus, we have taken care of both cases in the basis step.

For the inductive step, assume that $\HH{i}=1,$ for all $i$ such
that $k+1\leq i\leq s-1.$ Now let $v$ be one of the $h_s$ nodes of
type $v_k.$ One can show that the semigroup condition for $v$ in
the direction of $\Delta_{A}(v)$ cannot be satisfied if $\HH{k}\neq 1.$
In both cases $h_s=1$ and $h_s>1,$ the proof is essentially the
same as that of the basis step, so we omit the details.
\end{proof}

Proposition \ref{steps1and2} and Lemma \ref{induction} together
imply the following
\begin{corollary}\label{hscor} Suppose $h_s>1.$ Then the semigroup conditions imply that $\HH{k}=1$
for $1\leq k\leq s-1.$
\end{corollary}

In section \ref{toptypes}, we will see that for the case $h_s=1,$ the semigroup conditions \textit{and} congruence conditions together imply that $\HH{k}=1$
for $1\leq k\leq s-1.$

%%%%%%%%%%%%%%%%%%%%%%%%%%%%%%%%%%%%%%%%%%%%%%%%%%%%%%%%%%
%%%%%%%%%%%%%%%%%%%%% Section 5 %%%%%%%%%%%%%%%%%%%%%%%%%%%%%%%%
%%%%%%%%%%%%%%%%%%%%%%%%%%%%%%%%%%%%%%%%%%%%%%%%%%%%%%%%%%

\section{Action of the discriminant group}\label{discgroup_sec}

In order to use Proposition \ref{6.8} to check the congruence
conditions for the resolution graph $\Gamma_{f,n}$, we must compute
$e_w\cdot e_w$ for all leaves $w.$ By Proposition \ref{disc2},
this amounts to computing the continued fraction expansions of the
strings from leaves to nodes. This is essentially done in Mendris
and N\'{e}methi's paper (\cite{MenNem}, proof of Prop. 3.5), but
we need a bit more detail than they included.

\subsection{Background} We begin with a summary of facts that we
need, which can be found in \cite{coverings}.  Let $a,~Q,$
and $P$ be strictly positive integers with $\gcd(a,Q,P)=1.$  Let
$\sing{X(a,Q,P)}$ be the isolated surface singularity lying over
the origin in the normalization of $(\{U^aV^{Q}=W^{P}\},0).$ Let
$\lambda$ be the unique integer such that $0\leq \lambda< P/(a,P)$
and
$$Q+\lambda\cdot \frac{a}{(a,P)}=m\cdot\frac{P}{(a,P)},$$ for some
positive integer $m.$  If $\lambda\neq 0,$ then let
$k_1,\ldots,k_t\geq 2$ be the integers in the continued fraction
expansion of $\frac{P/(a,P)}{\lambda}.$ 

\begin{figure}[h]\centering\includegraphics[width=2.5in]{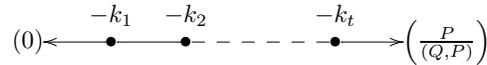}
\caption{The embedded resolution graph
$\Gamma(X(a,Q,P),V).$}\label{string}
\end{figure}

The minimal embedded resolution graph 
of the germ induced by the coordinate function $V$ on $\sing{X(a,Q,P)}$ is given by
the string in Figure \ref{string} (omitting the multiplicities of
the vertices).
If $\lambda=0,$ the string is empty.  One can
similarly describe the embedded resolution graphs of the functions
$U$ and $W,$ but we do not need them here.

\begin{lemma}\label{fractions} Let $N,~M,~P,$ and $Q$ be positive integers such that $(Q,P)=1$ and
$(N,M)=1.$ Let $\Gamma$ be the resolution graph of
the singularity in the normalization of $\sing{\{UV^Q=W^{P},~T^N=V^M\}}\subseteq\sing{\C^4}.$ Let $\lambda$ be the unique integer such that $0\leq
\lambda < P/(N,P)$ and
$$Q\frac{N}{(N,P)}+\lambda =m\cdot\frac{P}{(N,P)}$$ for some positive integer $m.$ Then if $\lambda\neq 0,$
$\Gamma$ is a string of vertices with continued fraction expansion
$\frac{P/(N,P)}{\lambda}.$  
\end{lemma}
\begin{proof}
We may assume $M=1,$ since it easy to check that the singularity in question has
the same normalization as
$\{UV^Q=W^{P},~T^N=V\}\subseteq\C^4.$
Therefore, $\Gamma$ is the resolution graph of the singularity in the
normalization of $\{UV^{QN}=W^{P}\},$ which is the same as the resolution graph of
$$X\left(1,Q\frac{N}{(N,P)},\frac{P}{(N,P)}\right)=\{UV^{QN/(N,P)}=W^{P/(N,P)}\}.$$\end{proof}

\subsection{Strings in $\Gamma_{f,n}$}\label{strings_sec} 
\begin{figure}[h]
\centering
\includegraphics[width=2in]{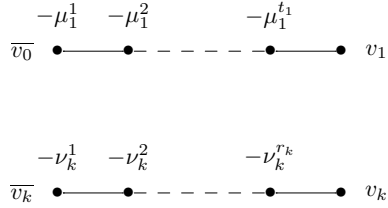}
\caption{Strings in
$\Gamma(\C^2,f)$.}\label{strings}
\end{figure}

We need the continued fraction expansion of
the strings in $\Gamma_{f,n}$ from leaves of type $\vbar{k},~0\leq k\leq
s,$ to the corresponding
node of type $v_k$ (from type $\vbar{0}$ to type $v_1$).  
First we recall the construction of $\Gamma(\C^2,f),$ the minimal good
embedded resolution graph of $f$ in $\C^2,$ as in
\cite{MenNem}. Let $f$ have Newton pairs $\{(p_k,q_k)~|~1\leq
k\leq s\}.$ Determine the continued fraction expansions
$$\frac{p_k}{q_k}=\mu_{k}^0-\dfrac{1}{\mu_{k}^1-\dfrac{1}{\ddots
-\dfrac{1}{\mu_{k}^{t_k}}}},\text{ and
}~\frac{q_k}{p_k}=\nu_{k}^0-\dfrac{1}{\nu_{k}^1-\dfrac{1}{\ddots
-\dfrac{1}{\nu_{k}^{r_k}}}},$$ where $\mu_{k}^0,~\nu_{k}^0\geq1,$
and $\mu_{k}^j,~\nu_{k}^j\geq2$ for $j>0.$ Then $\Gamma(\C^2,f)$ has the schematic form given in Figure \ref{schembres}. The strings from $\vbar{0}$ to $v_1$ and from $\vbar{k}$ to $v_k,~1\leq k\leq s,$ are given in Figure \ref{strings}. 
The multiplicities of the vertices $v_k$ are $m_{v_k}=a_k p_k p_{k+1}\cdots p_s, $ for $1\leq k\leq s.$

\begin{figure}[h]
\centering
\includegraphics[width=2.5in]{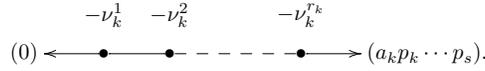}\caption{String from $\Gamma(\C^2,f)$.}\label{newstring}
\end{figure}

Consider the string in Figure \ref{newstring}.
The continued fraction expansion $[\nu^1_k,\ldots,\nu^{r_k}_k]$ corresponds to 
$p_k/\eta_k,$ where 
$q_k+\eta_k=\nu_k^0p_k.$ Let $X:=X(1,q_k,p_k).$ Then this string
is the embedded resolution graph of $V^{a_kp_{k+1}\cdots p_s}$ in
$X.$ It follows from the construction of $\Gamma_{f,n}$ that the
collection of strings that lies above this one in $\Gamma_{f,n}$ is the
(possibly non-connected) 
resolution graph of the singularity in the normalization of
$\{UV^{q_k}=W^{p_k},~T^n=V^{a_kp_{k+1}\cdots
p_s}\}.$  There are $(n,a_kp_{k+1}\cdots
p_s)=\h{k}d_k=\h{k}h_{k+1}\cdots h_s$ connected components (see
Definition \ref{allthestuff}), each being the resolution graph of
the normalization of
$$\{UV^{q_k}=W^{p_k},~T^{n/(\h{k}d_k)}=V^{a'_kp'_{k+1}\cdots
p'_s}\}.$$  Now we are in the situation of Lemma \ref{fractions},
with $Q=q_k,~P=p_k,$ and $N=n/(\h{k}d_k).$  We have
$(N,P)=(n/(\h{k}d_k),p_k)=h_k$ by definition of $h_k,$ and so in
this case $P/(N,P)=p'_k$ (as expected from Proposition
\ref{detmaxstrings}).  If $p'_k=1,$ then upon minimalization, the string
of type $\vbar{k}$ would completely collapse.

Suppose $p'_k\neq 1.$ By Lemma \ref{fractions}, the continued
fraction expansion of the string(s) from a leaf of type $\vbar{k}$
to the corresponding node of type $v_k$ in the minimalization of the resolution graph
$\Gamma_{f,n}$ is given by $p'_k/\eta'_k,$ where $\eta'_k$ is the unique
integer such that $0< \eta'_k <p'_k$ and
$$q_k \frac{n}{\HH{k}d_k}+\eta'_k=m p'_k,$$ for some positive integer $m.$
 Since $a_k=q_k+a_{k-1}p_{k-1}p_k,$ we have
\begin{equation}\label{etak'}\eta'_k\equiv -a_k \cdot \frac{n}{\HH{k}d_k}~(\text{mod
}p'_k).\end{equation} Knowing the congruence class of $\eta'_k$
modulo $p'_k$ is enough for our purposes.

The continued fraction expansion from $\vbar{0}$ to $v_1$ in $\Gamma(\C^2,f)$ is given
by $q_1/\eta_0=a_1/\eta_0,$ where 
$p_1+\eta_0=\mu_1^0 a_1.$ Using an argument analogous to the one
above, we have that if $a'_1\neq 1,$ the continued fraction expansion of the
string(s) from a leaf of type $\vbar{0}$ to the
corresponding node of type $v_1$ in the minimalization of $\Gamma_{f,n}$ is 
$a'_1/\eta'_0,$ where 
\begin{equation*}\label{eta0'}\eta'_0\equiv -p_1 \cdot
\frac{n}{\HH{1}d_1}~(\text{mod }a'_1).\end{equation*}

Recall the notation defined in section \ref{n-w}: for $r\in\Q,$
$[r]=\exp(2\pi i r),$ and for a leaf $w\in\Gamma_{f,n},$ $e_{w}$ denotes
the image in the discriminant group of the dual basis element in
$\mathbb{E}^*$ corresponding to $w.$

\begin{corollary}\label{edote} Let $w_k$ be any leaf of type $\vbar{k}$ in
$\Gamma_{f,n},~0\leq k\leq s,$ and assume that $p'_k\neq 1$ (assume $a'_1\neq 1$ for $k=0$).  Then
$$\left[e_{w_k}\cdot
e_{w_k}\right]=\left\{\begin{array}{ll}
\displaystyle\left[\frac{(n/\HH{1}d_1)(p_1a_2\cdots
a_s-A_1p'_1)}{a'_1a_2\cdots a_s}\right]&\text{ for }k=0\\
\displaystyle\left[\frac{(n/\HH{k}d_k)(a_ka_{k+1}\cdots
a_s-A_ka'_k)}{p'_ka_{k+1}\cdots a_s}\right]&\text{ for }1\leq
k\leq
s-1\\
\displaystyle\left[\frac{(n/\HH{s})(a_s-a'_s)}{p'_s}\right]&\text{
for }k=s.
\end{array}\right.$$
\end{corollary}

\begin{proof}  Proposition \ref{disc2} says that for a leaf $w$
connected by a string of vertices to a node $v,$ $$e_w\cdot
e_w=-d_v/(d^2\det(\Gamma))-p/d,$$ where $d_v$ is the product of
weights at the node $v,$ and $d/p$ is the fraction corresponding
to the string from $w$ to $v.$  Let $d_{v_k}$ be the product of
the weights at any node of type $v_k,~1\leq k\leq s$ (refer to
Figure \ref{sdk}). Then
$d_{v_k}=\da{k}\dm{k}^{h_k}(p'_k)^{\h{k}}.$ 

We need the following fact, which follows from Lemmas \ref{dm}
and \ref{det}. For any $k$ such that $1\leq k\leq s,$
$$\det(\Gamma_{f,n})=\frac{\dm{k}}{a'_k}\prod_{j=k}^s(p'_j)^{\h{j}-1}
\dm{j}^{h_j-1}.$$  
Now, for $1\leq k\leq s-1,$
\begin{eqnarray*}e_{w_k}\cdot e_{w_k}&=&-\frac{\da{k}\dm{k}^{h_k}(p'_k)^{\h{k}}}
{(p'_k)^2\det(\Gamma)}-\frac{\eta'_k}{p'_k}\\
&=&-\frac{ \left(\frac{n A_k \cdot \prod_{j=k+1}^{s}
(p'_j)^{\h{j}-1}\dm{j}^{h_j-1}}{\HH{k}d_ka_{k+1}\cdots
a_s}\right)\dm{k}^{h_k}(p'_k)^{\h{k}}}
{(p'_k)^2\frac{\dm{k}}{a'_k}\left(\prod_{j=k}^s(p'_j)^{\h{j}-1}
\dm{j}^{h_j-1}\right)}-\frac{\eta'_k}{p'_k}\\
&=&-\frac{n/(\HH{k}d_k)A_ka'_k} {p'_ka_{k+1}\cdots
a_s}-\frac{\eta'_k}{p'_k}.
\end{eqnarray*} Applying the congruence (\ref{etak'}), we have
$$[e_{w_k}\cdot e_{w_k}]=\left[\frac{(n/\HH{k}d_k)a_k}{p'_k}-
\frac{(n/\HH{k}d_k)A_ka'_k}{p'_ka_{k+1}\cdots
a_s}\right],
$$ and from here it is clear that the corollary is true. In the same way, it
is easy to check that that $\left[e_{w_0}\cdot e_{w_0}\right]$ and
$\left[e_{w_s}\cdot e_{w_s}\right]$ are as stated.\end{proof}

%%%%%%%%%%%%%%%%%%%%%%%%%%%%%%%%%%%%%%%%%%%%%%%%%%%%%%%%%%
%%%%%%%%%%%%%%%%%%%%% Section 6 %%%%%%%%%%%%%%%%%%%%%%%%%%%%%%%%
%%%%%%%%%%%%%%%%%%%%%%%%%%%%%%%%%%%%%%%%%%%%%%%%%%%%%%%%%%

\section{Proof of the Main Theorem}\label{toptypes}

In this section, we prove the Main Theorem, which determines
precisely which $\xfn,$ with $f$ irreducible, have a resolution
graph $\Gamma_{f,n}$ and associated splice diagram $\Delta_{f,n}$ that satisfy both
the
semigroup and congruence conditions.  

\begin{figure}[h]
\centering
\includegraphics[width=3.5in]{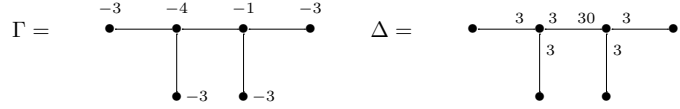}
\caption{Example for which the semigroup conditions are satisfied
but the congruence conditions are not.}\label{fig_a}
\end{figure}

\begin{remark}\label{mainrmk} \begin{itemize}\item[1)]The link is a $\Z$HS if and only if $n$ is relatively prime to all $p_i$ and $a_i$ (see \cite{zhs}). This is equivalent to all $h_i$ and $\h{i}$ being equal to $1.$ Hence this case belongs to (i)
of the Main Theorem. \item[2)]For the so-called pathological case
$n=p_s=2,$ both semigroup and congruence conditions are satisfied
only for $s=2.$  \item[3)]There are classes of $\xfn$ for which
the semigroup conditions are satisfied but the congruence
conditions are not, but we do not write up a complete list of
these types. An example with this property is given by
$n=2,~s=2,~p_1=2,~a_1=3,~p_2=3,$ and $a_2=20.$ The minimal good
resolution graph and splice diagram for this example are given in
Figure \ref{fig_a}.
\end{itemize}\end{remark}

We must treat
the cases $h_s=1$ and $h_s>1$ separately. The second case takes
much more work than the first.

\subsection{Case (i) $h_s=(n,p_s)= 1$}\label{toptypesi} 
First of all, we have the following
\begin{proposition}\label{hhk1} Suppose $h_s=1.$ If $\Gamma_{f,n}$ satisfies the semigroup and congruence conditions,
then $\HH{i}=1$ for $1\leq i\leq s-1.$
\end{proposition}

\begin{proof} 
In light of Lemma \ref{induction}, it suffices to
show that the semigroup and congruence conditions imply
$\HH{s-1}=1.$
We claim that the congruence condition at the unique node $v$ of
type $v_{s-1}$ cannot be satisfied if $\HH{s-1}\neq 1.$  Let
$u_j,~1\leq j\leq \h{s},$ denote the leaves of type $\vbar{s}$ in
$\Delta_{f,n},$ and let $y$ denote the leaf that arises from the string
$\Gamma(v_s)$ in $\Gamma^{can}(X_{f,n},z),$ as in Figure \ref{fig_d}.
\begin{figure}[h]
\centering
\includegraphics[width=3.5in]{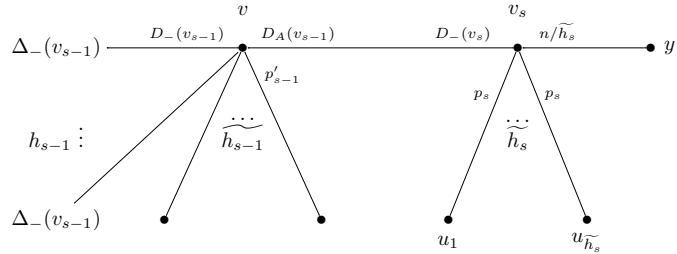}
\caption{Splice diagram for $h_s=1.$}\label{fig_d}
\end{figure} 
If $n/\h{s}=1,$ then the leaf $y$ does not exist,
but one can see that the argument holds regardless.

The semigroup condition at $v$ in the direction of $\Delta_A(v)$
says that there exist $\beta$ and $\alpha_i,~1\leq i\leq \h{s},$
in $\N\cup \{0\}$ such that
$$\da{s-1}=\left(\sum_{i=1}^{\h{s}}\alpha_i\right)(p_s)^{\h{s}-1}n/\h{s}
+\beta(p_s)^{\h{s}}.$$ It follows from Lemma \ref{dalemma} that 
$\da{s-1}=n/(\HH{s-1})(p_s)^{\h{s}-1}.$ Therefore, we have
\begin{equation}\label{admiss}n/(\HH{s-1})=\left(\sum_{i=1}^{\h{s}}\alpha_i\right)n/\h{s}+\beta
p_s.\end{equation} If $\h{s}=1,$ it is clear that $\HH{s-1}$ must
be $1;$ for, if not, $\alpha_1$ must be zero, which would imply
that $p_s$ divides $n/(\HH{s-1}).$ But this contradicts the
assumption that $h_s=1.$ Furthermore, note that if all $\alpha_i\geq
1,$ this implies that all $\alpha_i$ must equal 1, $\beta$ must be 0,
and $\HH{s-1}=1.$ If we assume $\HH{s-1}\neq 1,$ then there exists
$j$ such that $\alpha_j=0.$

Let $U_j$ be the variable associated to the leaf $u_j$
(respectively, $Y$ associated to $y$). By Proposition \ref{6.8},
the congruence condition at $v$ in the direction of $\Delta_A(v)$
implies, in particular, that there exists an admissible monomial
$H=U_1^{\alpha_1}\cdots U_{\h{s}}^{\alpha_{\h{s}}}Y^{\beta}$ such
that for every leaf $u_j,~1\leq j\leq \h{s},$
\begin{equation*}\left[\beta \frac{\ell_{yu_j}}{\det(\Gamma_{f,n})}+\sum_{i\neq
j}\alpha_i \frac{\ell_{u_iu_j}}{\det(\Gamma_{f,n})}- \alpha_j e_{u_j}\cdot
e_{u_j}\right]=\left[\frac{\ell_{vu_j}}{\det(\Gamma_{f,n})}\right].\end{equation*}
For the particular $j$ such that $\alpha_j=0,$ this condition is
\begin{equation}\label{step1cc}\left[\beta
\frac{\ell_{yu_j}}{\det(\Gamma_{f,n})}+\sum_{i\neq j}\alpha_i
\frac{\ell_{u_iu_j}}{\det(\Gamma_{f,n})}\right]=\left[\frac{\ell_{vu_j}}{\det(\Gamma_{f,n})}\right].\end{equation}
By Lemmas \ref{dm} and \ref{det},
$$\det(\Gamma_{f,n})=(p_s)^{\h{s}-1}\left(\frac{\dm{s}}{a'_{s}}\right)
=(p_s)^{\h{s}-1}(p'_{s-1})^{\h{s-1}-1}
\frac{\dm{s-1}^{h_{s-1}}}{a'_{s-1}}.$$ One can easily see that
$\left[\ell_{vu_j}/\det(\Gamma_{f,n})\right]=\left[0\right],$
$\left[\ell_{yu_j}/\det(\Gamma_{f,n})\right]=\left[0\right],$ and
$
\left[\ell_{u_iu_j}/\det(\Gamma_{f,n})\right]=\left[(a'_{s}n/\h{s})/p_s\right]~\text{for
}i\neq j.$ 
Thus the congruence condition (\ref{step1cc}) for the leaf $u_j$
is $\left[\left(\sum_{i\neq
j}\alpha_i\right)\frac{a'_{s}n/\h{s}}{p_s}\right]=\left[0\right];$
that is, $\left(\sum_{i\neq j}\alpha_i\right)a'_{s}n/\h{s}\in\Z p_s.$
Since $a'_s$ and $n/\h{s}$ are relatively prime to $p_s,$ this
implies that $\sum_{i\neq j}\alpha_i\in\Z p_s.$  But, by Equation
(\ref{admiss}), this implies that $n/(\HH{s-1})$ is divisible by
$p_s,$ which is a contradiction. Therefore, we must have $\HH{s-1}=1.$
\end{proof}

This leads us to the following 
\begin{proposition}\label{propi} Suppose $h_s=1.$  Then $\Gamma_{f,n}$ satisfies
the semigroup and congruence conditions if and only if both of the
following hold:\begin{itemize}\item[(I)]$\HH{i}=1$ for $1\leq i\leq
s-1,$
\item[(II)]$a'_s=a_s/\h{s}\in\N\langle a_{s-1},~p_1\cdots
p_{s-1},~a_{j}p_{j+1}\cdots p_{s-1}~:~1\leq j\leq
s-2\rangle.$\end{itemize}
\end{proposition}

\begin{remark}The condition (II) is clearly not always satisfied.  For example,
take $n$ divisible by $a_s.$ \end{remark}

\begin{proof} We have already shown (Propositions \ref{hhk1} and \ref{sgcleft})
that if the semigroup and congruence conditions are satisfied,
then (I) and (II) must hold. So assume that (I) and (II) are
satisfied. In the case that $\h{s}=1,$ the link is a $\zhs,$ and
the semigroup conditions are satisfied \cite{zhs}. (There are no congruence conditions when the link is a $\zhs.$)

Assume $\h{s}\neq 1.$ By Lemma \ref{dm}, $\dm{k}=a_k,~2\leq k\leq
s-1,$ and $ \dm{s}=a'_s,$ and it follows from Lemma \ref{dalemma} that
$\da{k}=n(p_s)^{\h{s}-1},$ for $1\leq k\leq s-1.$ There is exactly
one node of type $v_k$ in $\Delta_{f,n}$ for $1\leq k\leq s,$ which we
simply denote $v_k.$ We denote the leaves $z_0,\ldots, z_{s-1},$
$u_1,\ldots,u_{\h{s}},$ and $y,$ as in Figure \ref{fig_c}.
\begin{figure}[h]
\centering
\includegraphics[width=\textwidth]{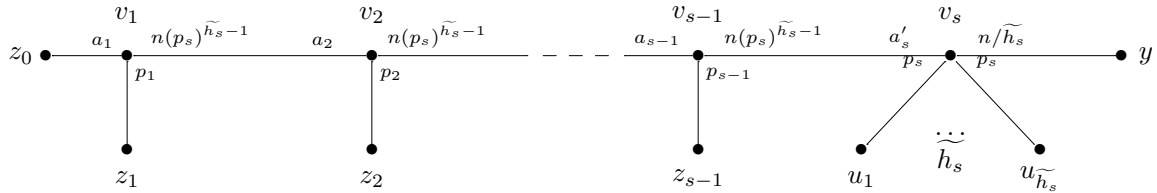}
\caption{Splice diagram for $\h{s}\neq 1$ and $\HH{i}=1,~1\leq
i\leq s-1.$}\label{fig_c}
\end{figure}

It is clear from Proposition \ref{sgcleft} that the semigroup
condition at the node $v_k$ in the direction of $\delm{k}$ is
satisfied for $2\leq k\leq s-1,$ and at the node $v_s,$ this
semigroup condition is equivalent to (II). Furthermore, one can
see by examination of the splice diagram that the semigroup
condition at each $v_k$ in the direction of $\dela{k}$ is always
satisfied (including in the case $n=\h{s}$).

It remains to show that $\Delta_{f,n}$ satisfies the congruence conditions.
Lemma \ref{det} implies that $\det(\Gamma_{f,n})=(p_s)^{\h{s}-1}.$  In
Figure \ref{fig_c}, it is easy to see that for any node $v$ and
any leaf $w$ in $\Delta_{f,n},$ $\ell_{vw}$ is always divisible by
$(p_s)^{\h{s}-1}.$ Therefore,
$\left[\ell_{vw}/\det(\Gamma_{f,n})\right]=[0]$ for any node $v$ and any
leaf $w.$ For each node, there are at most two conditions to
check: one for each adjacent edge that does not lead directly to a
leaf. By Proposition \ref{6.8}, we must show that for every node $v$ and
adjacent edge $e,$ there is an admissible
monomial $M_{ve}=\prod_{w\in\Delta_{ve}} Z_w^{\alpha_w}$ such that
for every leaf $w'$ in $\Delta_{ve},$
\begin{equation}\label{6.8eqn}\left[\sum_{w\neq
w'}\alpha_w \frac{\ell_{ww'}}{\det(\Gamma)}
-\alpha_{w'}e_{w'}\cdot
e_{w'}\right]=\left[0\right].\end{equation}

In this case, we have $A_i=a_{i+1}\cdots a_s$ for $1\leq i\leq
s-1.$ Since $A_1p'_1=a_2\cdots a_sp_1$ and $A_ja'_j=a_{j+1}\cdots
a_sa_j,$ Corollary \ref{edote} says that $[e_{z_j}\cdot
e_{z_j}]=[0]$ for $0\leq j\leq s-1.$ For any leaf $z_j,~0\leq
j\leq s-1,$ it is easy to see that $\ell_{z_jw'}$ is divisible by
$(p_s)^{\h{s}-1}$ for all leaves $w'\neq z_j$ in $\Delta_{f,n}.$ 
Since the subgraph $\delm{k}$ contains leaves only of the form
$z_j,~0\leq j\leq k-1,$ Equation (\ref{6.8eqn}) holds for all
leaves in $\delm{k}$ for any choice of admissible monomial. (In
fact, we have shown that the action of the discriminant group
element $e_{z_j}$ is trivial for $0\leq j\leq s-1.$)

Let $Z_{j}$ be the variable associated to the leaf $z_{j},$ $0\leq
j\leq s-1.$ It is easy to check that for $1\leq k\leq s-2,$ the
congruence condition at $v_k$ in the direction of $\dela{k}$ is
satisfied for
the admissible monomial $Z_{k+1}.$ 
The only remaining condition is for the node $v_{s-1}$ in the
direction of $v_{s}.$ Let $U_j$ be the variable associated to the
leaf $u_j,~1\leq j\leq \h{s}.$ We claim that the monomial
$U_1\cdots U_{\h{s}}$ (which is easily seen to be an admissible
monomial) satisfies the congruence condition.  It is clear from
the splice diagram that
$\left[\ell_{u_iu_j}/\det(\Gamma_{f,n})\right]=\left[(n/\h{s})a'_s/p_s\right]\text{
for }i\neq j,$ and by Corollary \ref{edote}, since each $u_j$ is a
leaf of type $\vbar{s},$ $[e_{u_j}\cdot
e_{u_j}]=\left[(n/\h{s})(a_s-a'_s)/p_s\right]\text{ for
 all }j.$  Hence, for each $u_j,$ Equation (\ref{6.8eqn})
 for the monomial $U_1\cdots U_{\h{s}}$ is
$$\left[(\h{s}-1)(n/\h{s})a'_s/p_s-(n/\h{s})(a_s-a'_s)/p_s\right]=[0].$$
This is clearly true, since $\h{s}a'_s=a_s.$  Finally, for the
leaf $y,$ Equation (\ref{6.8eqn})
 for $U_1\cdots U_{\h{s}}$ is
 $\left[\frac{\h{s}\ell_{yu_j}}{\det(\Gamma_{f,n})}
 \right]=\left[0\right]$ (for any choice of $j$). Since $\ell_{yu_j}$ is divisible by
$(p_s)^{\h{s}-1},$ the condition is
satisfied.\end{proof}

%%%%%%%%%%%%%%%%%%%%%%%%%%%%%%%%%%%%%%%%%%%%%%%%%%%%%%%%%%%%%%%%%%%%%%%%%%%%%%
%%%%%%%%%%%%%%%%%%%%%%%%%%%%%%%%%%%%%%%%%%%%%%%%%%%%%%%%%%%%%%%%%%%%%%%%%%%%%%
\subsection{Case (ii) $h_s=(n,p_s)> 1$}\label{toptypesii} The pathological case $n=p_s=2$
is treated separately at the end of the section. The main goal of
this section is to prove the following
\begin{proposition}\label{propii} Suppose $h_s>1$ and $n> 2.$ Then $\Gamma_{f,n}$ satisfies
the semigroup and congruence conditions if and only if
$$(\ast)~s=2,~p_2=2,~(n,p_2)=2,\text{ and }
(n,a_2)=(n/2,p_1)=(n/2,a_1)=1.$$
\end{proposition}

Let us first assume that $\Gamma_{f,n}$ satisfies the semigroup and
congruence conditions. We have already shown in $\S$\ref{sgc_sec}
that the semigroup conditions imply $h_s=(n,p_s)=p_s$ and
$\HH{i}=1$ for $1\leq i\leq s-1.$ Recall that since the link is a
$\Q$HS, $\h{s}=1$ and $a'_s=a_s.$  We prove that
$(\ast)$ must hold in two steps:\begin{itemize}\item[Step 1.]The
congruence conditions imply that $p_s=2.$ \item[Step 2.]The
congruence conditions imply that $s=2.$\end{itemize}

\begin{proof}[Proof of Step 1] For maximum convenience, we will use
the splice diagram $\Delta$ associated to the \textit{minimal}
good resolution graph $\Gamma^{min}(X_{f,n})$ (see Figure \ref{fig_g}). Recall
that $p'_s=1$ implies that there is no leaf of type $\vbar{s},$
since that string completely collapses in the minimal resolution
graph. We show that the congruence condition as in Proposition
\ref{6.8} for a node $v$ of type $v_{s-1}$ in the direction of
$\Delta_{A}(v)$ cannot hold unless $p_s=2.$ The only difficulty is in
notation.

By Lemmas \ref{dm} and \ref{dalemma}, $\dm{k}=a_k,\text{ for
}2\leq k\leq s,$ and
\begin{equation}\label{da6.2}\da{k}=\frac{n}{p_s}\tilde{A}_k(a_s)^{p_s-2},\text{ for
}1\leq k\leq s-1,\end{equation} where
$\tilde{A}_{s-1}=a_s-a_{s-1}p_{s-1}(p_s-1),$ and
$\tilde{A}_k=a_s-a_{k}p_{k}p_{k+1}^2\cdots p_{s-1}^2(p_s-1),$
$1\leq k\leq s-2.$
\begin{figure}[h]
\centering
\includegraphics[width=\textwidth]{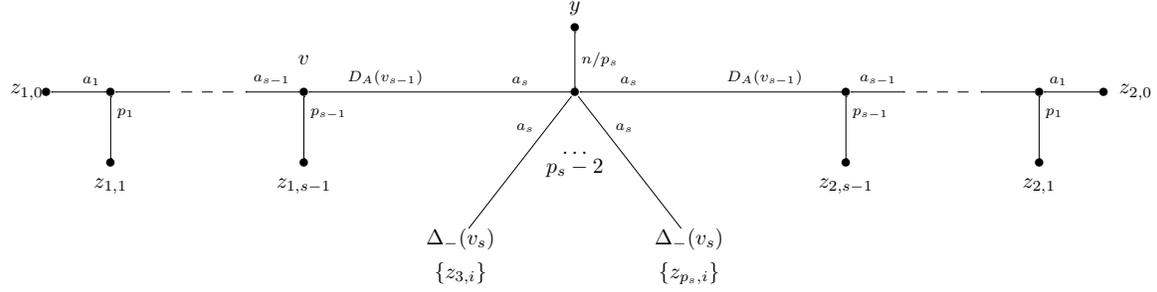}
\caption{Splice diagram for $h_s=p_s$ and $\HH{i}=1$ for $1\leq
i\leq s-1$.}\label{fig_g}
\end{figure}
Suppose that $p_s>2.$ For each $i,~0\leq i\leq s-1,$ there are
$h_s=p_s$ leaves of type $\vbar{i}.$  We label these leaves
$\{\z{j}{i}~|~1\leq j\leq p_s\},$ as indicated in Figure
\ref{fig_g}. The leaf on the edge with weight $n/p_s$ is denoted
$y,$ and is absent if $n/p_s=1.$ Let the corresponding variables
as in the Neumann-Wahl algorithm be $\{Z_{j,i}\}$ and $Y,$
respectively. Let $G$ be an admissible monomial for $v$ in the
direction of $\Delta_{A}(v)$ (i.e., in the direction of the
central node). We know that the variable $Y$ cannot appear in any
admissible monomial $G,$ by the proof of Proposition
\ref{steps1and2} ($M=0$). Therefore, we have
$G=\prod_{j=2}^{p_s}(Z_{j,0})^{\alpha_{j,0}}\cdots(Z_{j,s-1})^{\alpha_{j,s-1}},$
with $\alpha_{j,k}\in\N\cup\{0\}$ such that
\begin{equation}\label{bigsgc}\da{s-1}=\sum_{k=0}^{s-1}\sum_{j=2}^{p_s}\ell'_{v\z{j}{k}}\alpha_{j,k}.
\end{equation} 

For convenience of notation, we define integers $M_i$ as follows:
$$M_i:=\left\{
\begin{array}{ll}
p_1\cdots p_{s-1}& \text{for }i=0\\
a_ip_{i+1}\cdots p_{s-1}& \text{for }1\leq i\leq s-2\\
a_{s-1}& \text{for }i=s-1.
\end{array}\right.
$$ (Note that $M_i=\bar{\beta}_i/p_s.$) Let $v_s$ denote the unique node of type $v_s$ (the central node).
By Lemma \ref{genclaim}, $\ell'_{v_s\z{j}{i}}=M_i$ for all $j.$
Therefore, $\ell_{v\z{j}{i}}=M_ia_{s-1}p_{s-1}(a_s)^{p_s-2}n/p_s,$
and $\ell'_{v\z{j}{i}}=M_i(a_s)^{p_s-2}n/p_s,$
 for $1\leq i\leq s-1.$  Applying Equation (\ref{da6.2}) and
cancelling $(a_s)^{p_s-2}n/p_s$ from both sides of Equation
(\ref{bigsgc}) yields
\begin{equation}\label{bigsgc2}\tilde{A}_{s-1}=\sum_{k=0}^{s-1}\sum_{j=2}^{p_s}M_k
\alpha_{j,k}.
\end{equation}

Consider the congruence condition in Proposition \ref{6.8} for the
node $v$ in the direction of $\Delta_{A}(v)$ for each of the leaves
$\z{2}{i},~0\leq i\leq s-1.$ By Lemma \ref{det},
$\det(\Gamma_{f,n})=(a_s)^{p_s-1}.$ For any admissible monomial $G,$ the
condition for $w'=\z{2}{i}$ is equivalent to
\begin{equation}\label{bigcc}\left[\sum_{k=0}^{s-1}
\sum_{j=3}^{p_s}\alpha_{j,k}\frac{\ell_{\z{j}{k}\z{2}{i}
}}{(a_s)^{p_s-1}} +\sum_{k\neq
i}\alpha_{2,k}\frac{\ell_{\z{2}{k}\z{2}{i}
}}{(a_s)^{p_s-1}}-\alpha_{2,i}e_{\z{2}{i} }\cdot e_{\z{2}{i} }
\right]
=\left[\frac{\ell_{v\z{2}{i}}}{(a_s)^{p_s-1}} \right].\end{equation}   
For $0\leq i\leq s-1,$
\begin{equation}\label{lvz2i}\frac{\ell_{v\z{2}{i}}}{(a_s)^{p_s-1}}=
\frac{(n/p_s)M_ia_{s-1}p_{s-1}}{a_s}.\end{equation} Furthermore,
for any $j\neq2$
and for $0\leq k,i\leq s-1,$
\begin{equation}\label{lz1kz2i}\frac{\ell_{\z{j}{k}\z{2}{i} }}{(a_s)^{p_s-1}}=\frac{(n/p_s)
M_iM_k}{a_s}.\end{equation}

\begin{claim}\label{bigclaim} Fix $i$ such that $0\leq i\leq s-1.$
Then
\begin{itemize}\item[(a)] $\displaystyle[e_{\z{2}{i} }\cdot
e_{\z{2}{i} }]=\left[\frac{(n/p_s)M_i^2(p_s-1)}{a_s} \right],$ and
\item[(b)] For $k\neq i,$
$\displaystyle\left[\frac{\ell_{\z{2}{k}\z{2}{i}
}}{(a_s)^{p_s-1}}\right]=\left[\frac{-(n/p_s)M_iM_k(p_s-1) }{a_s}
\right],~0\leq k\leq s-1.$  \end{itemize}
\end{claim}

Let us assume for now that Claim \ref{bigclaim} is true and finish
the proof of Step 1.  By Equation (\ref{lz1kz2i}) and the Claim,
we have the following:
\begin{eqnarray*}
\text{Left side of
(\ref{bigcc})}&=&\left[\sum_{k=0}^{s-1}\sum_{j=3}^{p_s}\alpha_{j,k}
\frac{\frac{n}{p_s}M_iM_k}{a_s}
-\sum_{k=0}^{s-1}\alpha_{2,k}\frac{\frac{n}{p_s}M_iM_k(p_s-1)}{a_s}
\right]\\
&=&
\left[\frac{(n/p_s)M_i}{a_s}\left\{\sum_{k=0}^{s-1}\sum_{j=2}^{p_s}\alpha_{j,k}M_k
-p_s\sum_{k=0}^{s-1}\alpha_{2,k}M_k\right\} \right]\\
&=& \left[\frac{(n/p_s)M_i}{a_s}\left\{\tilde{A}_{s-1}
-p_s\sum_{k=0}^{s-1}\alpha_{2,k}M_k\right\}
\right]~(\text{by (\ref{bigsgc2}))}\\
&=& \left[\frac{(n/p_s)M_i}{a_s}\left\{a_s-a_{s-1}p_{s-1}(p_s-1)
-p_s\sum_{k=0}^{s-1}\alpha_{2,k}M_k\right\}
\right]\\
&=& \left[\frac{(n/p_s)M_i}{a_s}\left\{-a_{s-1}p_{s-1}(p_s-1)
-p_s\sum_{k=0}^{s-1}\alpha_{2,k}M_k\right\}
\right].\end{eqnarray*}

 Therefore, by (\ref{lvz2i}), the congruence condition
(\ref{bigcc}) is equivalent to
$$\left[\frac{\frac{n}{p_s}M_i}{a_s}\left\{-a_{s-1}p_{s-1}(p_s-1)
-p_s\sum_{k=0}^{s-1}\alpha_{2,k}M_k\right\}
\right]=\left[\frac{\frac{n}{p_s}M_ia_{s-1}p_{s-1}}{a_s}\right],$$
which is clearly equivalent to
$\left[-\frac{(n/p_s)M_ip_s}{a_s}\left(a_{s-1}p_{s-1}+
\sum_{k=0}^{s-1}\alpha_{2,k}M_k\right) \right]=[0].$ Since
$(a_s,n)=1$ and $(a_s,p_s)=1,$ this is equivalent to
\begin{equation}\label{finalcc}M_i\left(a_{s-1}p_{s-1}+
\sum_{k=0}^{s-1}\alpha_{2,k}M_k\right)\in \Z a_s.
\end{equation} Therefore, if the congruence conditions are
satisfied, that implies, in particular, that (\ref{finalcc}) holds
for all $i$ such that $0\leq i\leq s-1.$

We claim that if (\ref{finalcc}) holds for all $i,$ this implies
that $a_s$ divides $$\mathrm{S}:=a_{s-1}p_{s-1}+
\sum_{k=0}^{s-1}\alpha_{2,k}M_k. $$  Let $a_s=q_1^{e_1}\cdots q_l^{e_l}$ be the prime power factorization of $a_s.$ 
Suppose there is some $j$ such that $q_j^{e_j}$ does not divide $\mathrm{S}.$   Then at least one power of $q_j$
must divide $M_i$ for $0\leq i\leq s-1.$ In particular, $q_j$
divides $M_{s-1}=a_{s-1},$ and since $(a_{s-1},p_{s-1})=1,$ this
implies that $q_j$ divides $a_{s-2},$ because
$M_{s-2}=a_{s-2}p_{s-1}.$  This, in turn, implies $q_j$ divides
$a_{s-3},$ and so forth, down to $a_1.$ But $M_0=p_1\cdots
p_{s-1},$ which cannot possibly be divisible by $q_j.$  We have a contradiction, and thus $a_s$ divides $S.$

Finally, we claim that for $p_s>2,$ it is impossible for $a_s$ to
divide $\mathrm{S}.$  Equation (\ref{bigsgc2}), which is
equivalent to $a_s-a_{s-1}p_{s-1}(p_s-1)=\sum_{k=0}^{s-1}
\sum_{j=2}^{p_s}\alpha_{j,k}M_k,$ implies that
$\sum_{k=0}^{s-1}\alpha_{2,k}M_k\leq a_s-a_{s-1}p_{s-1}(p_s-1),$
and hence
$$\mathrm{S}=a_{s-1}p_{s-1}+\sum_{k=0}^{s-1}\alpha_{2,k}M_k\leq a_s-a_{s-1}p_{s-1}(p_s-2).$$
If $p_s>2,~a_s-a_{s-1}p_{s-1}(p_s-2)<a_s,$ which implies that
$\mathrm{S}<a_s,$ and hence $\mathrm{S}$ cannot be divisible by
$a_s,$ which is a contradiction. Therefore, we must have $p_s=2$
for the congruence conditions to be satisfied.\end{proof}

\begin{proof}[Proof of Claim \ref{bigclaim}] Since $\z{2}{i}$ is a leaf of type $\vbar{i},$ (a)
follows from Corollary \ref{edote}.
For (b), without loss of generality, we can assume $i<k.$ For
$1\leq i<k\leq s-2,~i\neq k-1,$ we have
$\ell_{\z{2}{k}\z{2}{i}}=\da{k}a_ip_{i+1}\cdots p_{k-1},$ and
hence,
\begin{eqnarray*}\left[\frac{\ell_{\z{2}{k}\z{2}{i}}}{\det(\Gamma_{f,n})} \right]
&=& \left[\frac{(n/p_s)\tilde{A}_ka_ip_{i+1}\cdots p_{k-1}}{a_s}
\right]\\
&=&\left[\frac{(n/p_s)(a_s-a_kp_kp_{k+1}^2\cdots p_{s-1}^2(p_s-1)
)a_ip_{i+1}\cdots p_{k-1}}{a_s} \right]\\
&=&\left[\frac{-(n/p_s)(p_s-1)a_kp_kp_{k+1}^2\cdots p_{s-1}^2
\cdot a_ip_{i+1}\cdots p_{k-1}}{a_s} \right]\\
&=&\left[\frac{-(n/p_s)(p_s-1)M_k M_i }{a_s} \right].
\end{eqnarray*}  The remaining cases 
are all similar and easy to check. \end{proof}

%%%%%%%%%%%%%%%%%%%%%%%%%%%%%%%%%%%%%%%%%%%%%%%%%%%%%%%%%%%%%%%%%%%%%%%%%%%%%%%%%%%%%%% Step 2
%%%%%%%%%%%%%%%%%%%%%%%%%%%%%%%%%%%%%%%%%%%%%%%%%%%%%%%%%%%%%%%%%%%%%%%%%%%%%%%%%%%%%%%

\begin{proof}[Proof of Step 2]So far, we have that the semigroup and congruence conditions
imply that $h_s=p_s=2$ and $\HH{i}=1$ for $1\leq i \leq s-1.$ Write
$n=2n'$ with $n'>1.$  We will show that for $s\geq 3,$ the
congruence conditions at a node $v$ of type $v_{s-2}$ in the
direction of $\Delta_{A}(v)$ cannot be satisfied. We should note that the
congruence condition at a node of type $v_{s-1}$ that we studied
in Step 1 \textit{can} be satisfied for $s\geq 3.$ For example,
take
$$\begin{array}{lll} a_1=3, & a_2=19, &
a_3=117,\\
p_1=2, & p_2=3, & p_3=2,
\end{array}$$ and any $n=2n'$ such that $n'$ is relatively prime to
$2,~3,~13,$ and $19.$

Figure \ref{fig_h} depicts the splice diagram in the general
situation.
\begin{figure}[h]
\centering
\includegraphics[width=\textwidth]{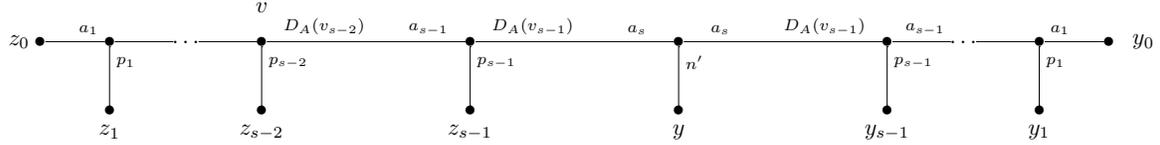}
\caption{Splice diagram for $n>2,~h_s=p_s=2,$ and $\HH{i}=1$ for
$1\leq i \leq s-1.$}\label{fig_h}
\end{figure}
The semigroup condition at $v$ in the direction of
$\Delta_{A}(v)$ is
$$\da{s-2}\in\N\langle \da{s-1},~a_sp_{s-1},~n'p_{s-1}M_i,~0\leq i\leq s-1\rangle.$$
Recall that $\da{s-1}=n'(a_s-a_{s-1}p_{s-1}),$ and
$\da{s-2}=n'(a_s-a_{s-2}p_{s-2}p_{s-1}^2).$ The semigroup condition implies that there exist 
$\alpha,~\beta,~\gamma_i\in\N\cup\{0\}$ such that 
\begin{equation*}n'(a_s-a_{s-2}p_{s-2}p_{s-1}^2)=\alpha
n'(a_s-a_{s-1}p_{s-1})+\beta a_sp_{s-1} +\sum_{i=0}^{s-1}\gamma_i
n'M_ip_{s-1}.\end{equation*} If $\beta\neq 0,$ then $\beta
a_sp_{s-1}$ must be divisible by $n'>1.$ By assumption,
$(a_s,n')=\h{s}=1,$ and $(p_{s-1},n')=h_{s-1}=1,$ and hence $n'$
must divide $\beta.$ But then $\beta a_sp_{s-1}\geq n'
a_sp_{s-1}>n'a_s>\da{s-2},$ and this is impossible.  Therefore,
$\beta=0.$

Hence, we can cancel $n'$ from the equation above, leaving
$$a_s-a_{s-2}p_{s-2}p_{s-1}^2=\alpha
(a_s-a_{s-1}p_{s-1})+\sum_{i=0}^{s-1}\gamma_i M_ip_{s-1}.$$  Since
$M_{s-1}=a_{s-1},$ we have
\begin{equation}\label{sgcblee}(\alpha-\gamma_{s-1})a_{s-1}p_{s-1}=(\alpha-1)a_s+\sum_{i=0}^{s-2}\gamma_i
M_ip_{s-1}+a_{s-2}p_{s-2}p_{s-1}^2,\end{equation} which implies
$(\alpha-\gamma_{s-1})a_{s-1}p_{s-1}>(\alpha-1)a_s.$  Suppose
$\alpha>1.$  Then, since $a_s=q_s+a_{s-1}p_{s-1}p_s$ and $p_s=2,$
$$(\alpha-\gamma_{s-1})a_{s-1}p_{s-1}>(\alpha-1)a_s>(\alpha-1)2a_{s-1}p_{s-1}.$$
  This implies
$(\alpha-\gamma_{s-1})-2(\alpha-1)>0,$ i.e.,
$2>\alpha+\gamma_{s-1}.$  But this is impossible for $\alpha>1.$

Now suppose $\alpha=1.$ It is clear from Equation (\ref{sgcblee})
that $\gamma_{s-1}$ must be $0,$ and so we have
$$a_{s-1}p_{s-1}=\sum_{i=0}^{s-2}\gamma_i
M_ip_{s-1}+a_{s-2}p_{s-2}p_{s-1}^2,$$ i.e.,
$a_{s-1}=\sum_{i=0}^{s-2}\gamma_i M_i+a_{s-2}p_{s-2}p_{s-1}.$  But
$M_i$ is divisible by $p_{s-1}$ for $0\leq i\leq s-2,$ so this
would imply $a_{s-1}$ is divisible by $p_{s-1},$ which is
impossible. Therefore, $\alpha=0,$ and we have
\begin{equation}\label{sgc4}a_s-a_{s-2}p_{s-2}p_{s-1}^2=\sum_{i=0}^{s-1}\gamma_i
M_ip_{s-1}.\end{equation}  (Note that this semigroup condition is
already quite restrictive, because it requires $a_s$ to be
divisible by $p_{s-1}.$)

Now let us return to the congruence conditions for the node $v$ in
the direction of $\Delta_{A}(v).$  An admissible monomial for $v$ in that
direction must be of the form $H=Y_0^{\gamma_0}\cdots
Y_{s-1}^{\gamma_{s-1}},$ with $\gamma_i\in\N\cup\{0\}.$ The
congruence condition for the leaf $y_{s-1}$ is
$$\left[\frac{\ell_{vy_{s-1}}}{\det(\Gamma_{f,n})}\right]=
\left[\sum_{i=0}^{s-2}\gamma_i\frac{\ell_{y_{s-1}y_i}}{\det(\Gamma_{f,n})}-\gamma_{s-1}e_{y_{s-1}}\cdot
e_{y_{s-1}}\right].$$  
Applying  Claim \ref{bigclaim}, this condition is equivalent to
$$\left[\frac{n'a_{s-2}p_{s-2}a_{s-1}p_{s-1}}{a_s}\right]=\left[-\frac{n'a_{s-1}}{a_s}
\left(\sum_{i=0}^{s-1}\gamma_iM_i\right)\right];$$ that is,
$n'a_{s-1}\left(a_{s-2}p_{s-2}p_{s-1}+\sum_{i=0}^{s-1}\gamma_iM_i\right)\in\Z
a_s.$  Since $(a_s,n')=1,$ we must have\\
$a_{s-1}\left(a_{s-2}p_{s-2}p_{s-1}+\sum_{i=0}^{s-1}\gamma_iM_i\right)=Na_s$
for some $N$ in $\Z.$  If we multiply both sides of this equation
by $p_{s-1}$ and apply Equation (\ref{sgc4}), we get
$$a_{s-1}a_{s-2}p_{s-2}p_{s-1}^2+a_{s-1}(a_s-a_{s-2}p_{s-2}p_{s-1}^2)=Na_sp_{s-1};$$ i.e., $a_{s-1}=Np_{s-1}.$  This
implies $p_{s-1}$ divides $a_{s-1},$ which is a contradiction.

Therefore, we have shown that if $s\geq 3,$ then the congruence
condition for the node $v$ of type $v_{s-2}$ in the direction of
$\Delta_{A}(v)$ cannot be satisfied for the leaf $y_{s-1}.$  Hence, the
congruence conditions imply that $s=2.$\end{proof} We have
finished Steps 1 and 2, hence have proved one direction of
Proposition \ref{propii}.

%%%%%%%%%%%%%%%%%%%%%%%%%%%%%%%%%%%%%%%%%%%%%%%%%%%%%%%%%%%%%%%%%%%%%%%%%%%%%%%%%%%%%%%%%%%%%%%
%%%%%%%%%%%%%%%%%%%%%%%%%%%%%%%%%%%%%%%%%%%%%%%%%%%%%%%%%%%%%%%%%%%%%%%%%%%%%%%%%%%%%%%%%%%%%%%

\begin{figure}[h]
\centering
\includegraphics[width=3in]{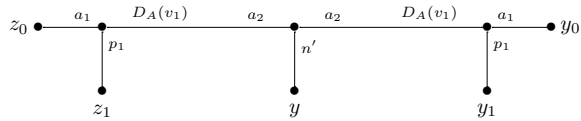}
\caption{Splice diagram for ($\ast$), $n>2$.}\label{fig_e}
\end{figure}

For the other direction, we must check that ($\ast$) implies that
the semigroup and congruence conditions are satisfied. The splice
diagram in this situation is shown in Figure \ref{fig_e}.
The only semigroup condition that needs to be checked is
$$\da{1}\in\N\langle a_2,~n'a_1,~n'p_1 \rangle,$$ where $\da{1}=n'(a_2-a_1p_1)=n'(q_2+a_1p_1).$
Since $a_1$ and $p_1$ are relatively prime, the conductor of the
semigroup generated by $a_1$ and $p_1$ is less than $a_1p_1,$
hence $a_1p_1+q_2$ is in the semigroup generated by $a_1$ and
$p_1,$ and therefore this semigroup condition is satisfied. 

There are only two congruence conditions to check. One is
equivalent to the following: there exist $\alpha_0$ and $\alpha_1$
in $\N\cup\{0\}$ such that $a_2=\alpha_0 p_1+\alpha_1 a_1,$
$$\left[\alpha_1\frac{-n'a_1p_1}{a_2}-\alpha_0\frac{n'p_1^2}{a_2}\right]=[0],\text{ and }
\left[\alpha_0\frac{-n'a_1p_1}{a_2}-\alpha_1\frac{n'a_1^2}{a_2}\right]=[0].$$
But these conditions are obviously both satisfied for any
$\alpha_0,~\alpha_1$ such that $a_2=\alpha_0 p_1+\alpha_1 a_1.$
The other congruence condition is equivalent to the following:
there exist $\gamma_0$ and $\gamma_1$ in $\N\cup\{0\}$ such that
$a_2-a_1p_1=\gamma_0 p_1+\gamma_1 a_1,$
$$\left[\gamma_1\frac{-n'a_1p_1}{a_2}-\gamma_0\frac{n'p_1^2}{a_2}\right]=
\left[\frac{n'a_1p_1^2}{a_2}\right],\text{ and }\left[\gamma_0\frac{-n'a_1p_1}{a_2}-\gamma_1\frac{n'a_1^2}{a_2}\right]=
\left[\frac{n'a_1^2p_1}{a_2}\right].$$ But these conditions are
also obviously both satisfied for any $\gamma_0,~\gamma_1$ such
that $a_2-a_1p_1=\gamma_0 p_1+\gamma_1 a_1.$ This concludes the
proof of Proposition \ref{propii}.

\subsubsection*{The pathological case} If $h_s>1$ and $n=2,$ then the semigroup conditions imply that $p_s=2$ by Proposition \ref{steps1and2}. Therefore, all that remains in the
proof of the Main Theorem is the pathological case. Let $\Gamma_{f,n}$ be
the graph associated to the minimal good resolution (see
$\S$\ref{top}).

\begin{proposition}Suppose $n=p_s=2.$ Then $\Gamma_{f,n}$ satisfies
the semigroup and congruence conditions if and only if
$s=2.$\end{proposition}

\begin{proof}We begin by assuming that $\Gamma_{f,n}$ satisfies the
semigroup and congruence conditions.  It is automatically true that $\HH{i}=1$ for $1\leq i\leq s-1,$ and
that $h_s=2.$  We must show that $s$ must be $2.$ The
splice diagram is pictured in Figure \ref{fig_i}.
\begin{figure}[h]
\centering
\includegraphics[width=\textwidth]{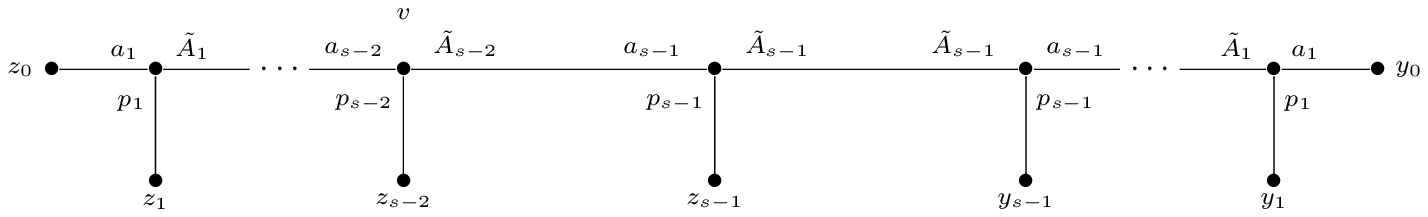}
\caption{Splice diagram for the pathological case,
$s>2.$}\label{fig_i}
\end{figure}
We can use essentially the same argument as in Step 2 above to
show that for $s\geq 3,$ the congruence conditions at the node $v$
of type $v_{s-2}$ in the direction of $\Delta_{A}(v)$ cannot possibly be
satisfied for the leaf $y_{s-1}.$

The semigroup condition at $v$ in the direction of $\Delta_{A}(v)$ is
$$\tilde{A}_{s-2}\in\N\langle \tilde{A}_{s-1},~p_{s-1}M_i,~0\leq i\leq s-1\rangle.$$
Precisely the same argument as in Step 2 above shows that
$\tilde{A}_{s-1}$ cannot appear in the expression for $\tilde{A}_{s-2}$ that
comes from the semigroup condition. Therefore, there exist
$\gamma_i$ in $\N\cup\{0\}$ such that
$a_s-a_{s-2}p_{s-2}p_{s-1}^2=\sum_{i=0}^{s-1}\gamma_i M_ip_{s-1}.$

Let $H=Y_0^{\gamma_0}\cdots Y_{s-1}^{\gamma_{s-1}}$ be an
admissible monomial for $v$ in the direction of $\Delta_{A}(v).$  The
congruence
condition for the leaf $y_{s-1}$ 
is equivalent to
$$\left[\frac{a_{s-2}p_{s-2}a_{s-1}p_{s-1}}{a_s}\right]=\left[-\frac{a_{s-1}}{a_s}
\left(\sum_{i=0}^{s-1}\gamma_iM_i\right)\right].$$ 
Just as in Step 2, this implies $p_{s-1}$ divides $a_{s-1},$ and
hence the congruence conditions cannot be satisfied for $s>2.$

Finally, for $s=2,$ it is easy to check that the semigroup and
congruence conditions are satisfied.\end{proof}

%%%%%%%%%%%%%%%%%%%%%%%%%%%%%%%%%%%%%%%%%%%%%%%%%%%%%%%%%%%%%%%
%%%%%%%%%%%%%%%%%%%%%%%%%%% Bibliography %%%%%%%%%%%%%%%%%%%%%%%%%%%%%
%%%%%%%%%%%%%%%%%%%%%%%%%%%%%%%%%%%%%%%%%%%%%%%%%%%%%%%%%%%%%%%

\bibliographystyle{amsplain}
\bibliography{biblio_mmj}

\end{document}